\documentclass[11pt,a4paper]{amsart}
\usepackage{amsmath}
\usepackage{graphicx}
\usepackage{epsfig}
\usepackage{xypic}

\newtheorem{prop}{Proposition}[section]
\newtheorem{lemma}{Lemma}[section]
\newtheorem{cor}{Corollary}[section]
\newtheorem{defi}{Definition}[section]
\newtheorem{theorem}{Theorem}[section]

\input epsf
\input cyracc.def

\newcommand{\bprf}{{\it Proof.~}}
\newcommand{\eprf}{\hfill $\square$ \bigskip\par}
\newcommand{\C}{\mathbb{C}}
\newcommand{\Q}{\mathbb{Q}}

\newcommand{\Z}{\mathbb{Z}}
\newcommand{\PP}{\mathbb{P}}
\newcommand{\N}{\mathbb{N}}
\newcommand{\cl}{\mathcal{L}}
\newcommand{\lra}{\longrightarrow}
\makeatletter                                                  
\def\blfootnote{\xdef\@thefnmark{}\@footnotetext}                   
\begin{document}
\title{Symplectic automorphisms of prime order on \\ K3 surfaces}
\author{Alice Garbagnati and Alessandra Sarti}
\address{Alice Garbagnati, Dipartimento di Matematica, Universit\`a di Milano,
  via Saldini 50, I-20133 Milano, Italy}

\email{garba@mat.unimi.it}

\address{Alessandra Sarti, Institut f\"ur  Mathematik, Universit\"at Mainz,
Staudingerweg 9, 55099 Mainz, Germany}

\email{sarti@mathematik.uni-mainz.de}

\begin{abstract}
We study algebraic K3 surfaces (defined
over the complex number field) with a symplectic automorphism of
prime order. In particular we consider the action of the
automorphism on the second cohomology with integer coefficients (by a result of Nikulin this action is independent on the choice of the K3 surface).
With the help of elliptic fibrations we determine the invariant sublattice and its perpendicular
complement, and show that the latter coincides with the
Coxeter-Todd lattice in the case of automorphism of order three.

\end{abstract}

\maketitle

\pagestyle{myheadings}

\markboth{ ALICE GARBAGNATI AND ALESSANDRA SARTI}{SYMPLECTIC AUTOMORPHISMS OF PRIME ORDER ON K3 SURFACES}

\blfootnote {The second author was partially supported by DFG Research Grant SA 1380/1-2.}

\section{Introduction}
In the paper \cite{Nikulin} Nikulin studies finite abelian groups
$G$ acting symplectically (i.e. $G_{|H^{2,0}(X,\mathbb{C})}=
id_{|H^{2,0}(X,\mathbb{C})}$) on K3 surfaces (defined over $\C$).
One of his main result is that the action induced by $G$ on the
cohomology group $H^2(X,\mathbb{Z})$ is unique up to isometry.
In \cite{Nikulin} all abelian finite groups of automorphisms of a
K3 surface acting symplectically are classified. Later Mukai in
\cite{Mukai} extends the study to the non abelian case. Here we
consider only abelian groups of prime order $p$
which, by Nikulin, are isomorphic to $\mathbb{Z}/p\mathbb{Z}$ for $p=2,3,5,7$.\\
In the case of $p=2$ the group is generated by an involution,
which is called by Morrison in \cite[Def. 5.1]{morrison} Nikulin
involution. This was very much studied in the last years, in
particular because of its relation with the Shioda-Inose structure
(cf. e.g. \cite{doran}, \cite{galluzzi},  \cite{geemen},
\cite{long}, \cite{morrison}). In \cite{morrison} Morrison proves
that the isometry induced by a Nikulin involution $\iota$ on the
lattice $\Lambda_{K3}\simeq U\oplus  U\oplus U\oplus E_8(-1)\oplus
E_8(-1)$, which is isometric to $H^2(X,\Z)$, switches the two
copies of $E_8(-1)$ and acts as the identity on the sublattice
$U\oplus U\oplus U$. As a consequence one sees that
$(H^2(X,\mathbb{Z})^{\iota^*})^{\perp}$ is the lattice $E_8(-2)$.
This implies that the Picard number $\rho$ of an algebraic K3
surface admitting a Nikulin involution is at least nine. In
\cite{bert} van Geemen and Sarti show that if $\rho\geq 9$ and
$E_8(-2)\subset NS(X)$ then the algebraic $K3$ surface $X$ admits
a Nikulin involution and they classify completely these $K3$
surfaces. Moreover they discuss  many examples and in particular
those surfaces admitting an elliptic fibration with a section of
order two.
This section operates by translation on the fibers and defines a Nikulin involution on the K3 surface.\\
The aim of this paper is to identify the action of a symplectic
automorphism $\sigma_p$ of the remaining possible prime orders
$p=3,5,7$ on the K3 lattice $\Lambda_{K3}$ and to describe such
algebraic K3 surfaces with minimal possible Picard number. Thanks
to Nikulin's result (\cite[Theorem 4.7]{Nikulin}), to find the
action on $\Lambda_{K3}$, it suffices to identify the action in
one special case.  For this purpose it seemed to be convenient to
study algebraic K3 surfaces with an elliptic fibration
 with a section of order three, five, resp. seven. Then
the translation by this section is a symplectic automorphism of
the surface of the same order. A concrete analysis leads us to the
main result of the paper which is the description of the lattices
$H^2(X,\mathbb{Z})^{\sigma_p^*}$ and
$\Omega_p=(H^2(X,\mathbb{Z})^{\sigma_p^*})^{\perp}$ given in the
Theorem \ref{main}. The proof of the main theorem consists in the
Propositions \ref{invariant and orthogonal part 3}, \ref{invariant
and orthogonal part 5}, \ref{invariant and orthogonal part 7}. We
describe the lattice $\Omega_p$ also as
$\mathbb{Z}[\omega_p]$-lattices, where $\omega_p$ is a primitive $p$ root of
the unity. This kind of lattices are studied e.g. in
\cite{bayer1}, \cite{bayer2} and \cite{Ebeling}. In particular in
the case $p=3$ the lattice $\Omega_3$ is the Coxeter-Todd lattice
with the form multiplied by $-2$, $K_{12}(-2)$, which is described
in
\cite{Coxeter Todd} and in \cite{Conway Sloane}.\\
The elliptic surfaces we used to find the lattices $\Omega_p$ do
not have the minimal possible Picard number. We prove in
Proposition \ref{latticeminimalrank} that for K3 surfaces, $X$,
with minimal Picard number and symplectic automorphism, if $L$ is
a class in $NS(X)$ which is invariant for the automorphisms,
with $L^2=2d>0$,  then either
$NS(X)=\mathbb{Z}L\oplus\Omega_p$ or  the latter is a sublattice
of index $p$ in $NS(X)$. Using this result and the one of the
Proposition \ref{periodi} we describe the coarse moduli space of
the algebraic K3 surfaces admitting a symplectic automorphism
of prime order.\\
The structure of the paper is the following: in section
\ref{section:moduli} we compute the number of moduli of algebraic
K3 surfaces admitting a symplectic automorphism of order $p$ and
their minimal Picard number. In section \ref{z[omega] lattices} we
give the definition of $\mathbb{Z}[\omega_p]$-lattice and we
associate to it a module with a bilinear form,
 which in some cases is a $\mathbb{Z}$-lattice
(we use this construction in section \ref{elliptic} to describe
the lattices $\Omega_p$ as $\mathbb{Z}[\omega_p]$-lattices). In
section \ref{elliptic: general facts} we recall some results about
elliptic fibrations and elliptic K3 surfaces (see e.g.
\cite{miranda1}, \cite{miranda2}, \cite{shimada}, \cite{shioda}
for more on elliptic K3 surfaces). In particular we introduce the
three elliptic fibrations which we use in section \ref{elliptic}
and give also their Weierstrass form. In section \ref{elliptic} we
state and proof the main result, Theorem \ref{main}: we identify
the lattices $\Omega_p$ and we describe them as
$\mathbb{Z}[\omega_p]$-lattices. In section \ref{minpic} we
describe the N\'eron-Severi group of K3 surfaces admitting a
symplectic automorphism and having minimal Picard number
(Proposition \ref{latticeminimalrank}). In section 5 we describe
the coarse moduli space of the algebraic K3 surfaces admitting a
symplectic automorphism and the N\'eron-Severi group of those
having minimal Picard number.\\
{\it We would like to express our deep thanks to Bert van Geemen
for suggesting us the problem and for his invaluable help during
the preparation of this paper.}

\section{Preliminary results}\label{section:moduli}
\begin{defi}
A {\it symplectic automorphism} $\sigma_p$ of order $p$ on a K3
surface $X$ is an automorphism such that:

\smallskip

1. the group $G$ generated by $\sigma_p$ is isomorphic to $\Z/p\Z$,\\
2. $\sigma_p^*(\delta)=\delta$, for all $\delta$ in $H^{2,0}(X).$
\end{defi}

We recall that by \cite{Nikulin} an automorphism on a K3 surface
is symplectic if and only if it acts as the identity on the
transcendental lattice $T_X$. In local coordinates at a fixed
point $\sigma_p$ has the form
$\mbox{diag}(\omega_p,\omega_p^{p-1})$ where $\omega_p$ is  a primitive
$p$-root of unity. By a result of Nikulin the only possible values
for $p$ are 2,3,5,7 see \cite[Theorem 4.5]{Nikulin} and \cite[\S
5]{Nikulin}. The automorphism $\sigma_3$ has six fixed points on
$X$, $\sigma_5$ has four fixed points and $\sigma_7$ has three
fixed points. The automorphism $\sigma_p$ induces a $\sigma_p^*$
isometry on $H^2(X,\Z)\cong \Lambda_{K3}$. Nikulin proved
\cite[Theorem 4.7]{Nikulin} that if $\sigma_p$ is symplectic, then
the action of
$\sigma_p^*$ is unique up to isometry of $\Lambda_{K3}$.\\
Let $\omega_p$ be a primitive
$p$-root of the unity. The vector space $H^2(X,\mathbb{C})$ can be
decomposed in eigenspaces of the eigenvalues 1 and $\omega_p^i$:
$$
H^2(X,\mathbb{C})=H^2(X,\mathbb{C})^{\sigma_p^*}\oplus(\bigoplus_{i=1,\ldots,p-1}H^2(X,\mathbb{C})_{\omega_p^i}).
$$
We observe that the non rational eigenvalues $\omega_p^i$ have all
the same multiplicity. So we put:
\smallskip
$a_p:=$ multiplicity of the eigenvalue 1,\ \  $b_p:=$multiplicity
of the eigenvalues $\omega_p^i$.\\
\smallskip
In the following we find $a_p$ and $b_p$ by using the Lefschetz fixed
point formula:
\begin{eqnarray}\label{Lefschetz fixed point}
\mu_p=\sum_r(-1)^r {\rm trace}(\sigma_p^*|H^r(X,\mathbb{C}))
\end{eqnarray}
where $\mu_p$ denotes the number of fixed points. For $K3$
surfaces we obtain
\begin{eqnarray*}
\mu_p=1+0+{\rm trace}(\sigma_p^*|H^2(X,\C))+0+1.
\end{eqnarray*}

\begin{prop}\label{theorem moduli}
Let $X$, $\sigma_p$, $a_p$, $b_p$ be as above, $p=3,5,7$. Let
$\rho_p$ be the Picard number of $X$, and let $m_p$ be the
dimension of the moduli space of the algebraic K3 surfaces
admitting a symplectic automorphism of order $p$. Then
$$
\begin{array}{llll}
a_3=10\ \ \ &b_3=6\ \ \ \ \ \ \ &\rho_3\geq13\ \ \ \ \ &m_3\leq 7\\
a_5=6\ \ \ &b_5=4 \ \ \ \ \ \ \ &\rho_5\geq17\ \ \ \ \ &m_5\leq 3\\
a_7=4\ \ \ &b_7=3\ \ \ \ \ \  \ &\rho_7\geq19\ \ \ \ \ &m_7 \leq 1.
\end{array}
$$
\end{prop}
\bprf The proof is similar in all the cases, here we give the details only in the case $p=5$.\\
A symplectic automorphism of order five on a K3 surface has
exactly four fixed points. Applying the Lefschetz fixed points
formula \eqref{Lefschetz fixed point}, we have
$a_5+b_5(\omega_5+\omega_5^2+\omega_5^3+\omega_5^4)=2$. Since
$\omega_p^{p-1}=-(\sum_{i=0}^{p-2}\omega_p^i)$, the equation
becomes $a_5-b_5=2$.\\
Since $\dim H^2(X,\mathbb{C})=22$, $a_5$ and $b_5$ have to satisfy:
\begin{eqnarray}\label{system a_p b_p}
\left\{\begin{array}{ccc}
a_5-b_5&=&2\\
a_5+4b_5&=&22.
\end{array}\right.
\end{eqnarray}
We have $$\begin{array}{l} \dim H^2(X,\C)^{\sigma_5^*}=6=a_5\ \mbox{and}\\
\dim H^2(X,\C)_{\omega_5}=\dim H^2(X,\C)_{\omega_5^2}=\dim
H^2(X,\C)_{\omega_5^3}=\dim H^2(X,\C)_{\omega_5^4}=4=b_5.
\end{array}$$
Since  $T_X\otimes\C \subset H^2(X,\C)^{\sigma_5^*}$,
$(H^2(X,\C)^{\sigma_5^*})^{\perp}=H^2(X,\C)_{\omega_5}\oplus
H^2(X,\C)_{\omega_5^2}\oplus H^2(X,\C)_{\omega_5^3}\oplus
H^2(X,\C)_{\omega_5^4}\subset NS(X)\otimes \C$. We consider only
algebraic $K3$ surfaces and so we have an ample class $h$ on X, by
taking $h+\sigma_5^* h+\sigma_5^{*2} h+\sigma_5^{*3}
h+\sigma_5^{*4} h$ we get a $\sigma_5$-invariant class, hence in
$H^2(X,\C)^{\sigma_5^*}$. From here it follows that $\rho_p=$rank
$NS(X) \geq 16+1=17$, whence rank $T_X\leq 22-17=5$. The number of
moduli is at most $20-17=3$.\eprf

\textbf{Remark.} In \cite[\S 10]{Nikulin} Nikulin computes
$\mbox{rank}(H^2(X,\mathbb{Z})^{\sigma_p^*})^{\perp}=(p-1)b_p$ and
$\mbox{rank}(H^2(X,\mathbb{Z})^{\sigma_p^*})=a_p$. In \cite[Lemma
4.2]{Nikulin} he also proves that there are no classes with self
intersection $-2$ in the lattices
$(H^2(X,\mathbb{Z})^{\sigma_p^*})^{\perp}$; we describe these
lattices in the sections \ref{elliptic3}, \ref{elliptic5},
\ref{elliptic7} and we find again the result of Nikulin.

\section{The $\Z[\omega]$-lattices}\label{z[omega] lattices}

In the sections \ref{section lattice omega_3}, \ref{section
lattice omega_5}, \ref{section lattice omega_7} our purpose is to
describe $(H^2(X,\mathbb{Z})^{\sigma_p^*})^{\perp}$ as
$\mathbb{Z}[\omega_p]$-lattice. We recall now some useful results
on these lattices.

\begin{defi}
Let $p$ be an odd prime and $\omega:=\omega_p$ be a  primitive $p$-root of
the unity. A $\mathbb{Z}[\omega]$-lattice is a free
$\mathbb{Z}[\omega]$-module with an hermitian form (with values
in $\mathbb{Z}[\omega]$). Its rank is its rank as
$\mathbb{Z}[\omega]$-module.
\end{defi}

Let $\{L,h_L\}$ be a $\mathbb{Z}[\omega]$-lattice of rank $n$. The
$\mathbb{Z}[\omega]$-module $L$ is also a $\mathbb{Z}$-module of
rank $(p-1)n$. In fact if $e_i$, $i=1,\ldots, n$ is a basis of $L$
as $\mathbb{Z}[\omega]$-module, $\omega^je_i$, $i=1,\ldots,n$,
$j=0,\ldots,p-2$ is a basis for $L$ as $\mathbb{Z}$-module (recall
that $\omega^{p-1}=-(\omega^{p-2}+\omega^{p-3}+\ldots+1)$). The
$\mathbb{Z}$-module $L$ will be called $L_{\mathbb{Z}}$.\\
Let $\Gamma_p:=Gal(\mathbb{Q}(\omega)/\mathbb{Q})$ be the group of the
automorphisms of $\mathbb{Q}(\omega)$ which fix $\mathbb{Q}$.  We
recall that the group $\Gamma_p$
has order $p-1$ and its elements are automorphisms $\rho_i$ such
that $\rho_i(1)=1,\ \ \rho_i(\omega)=\omega^i$ where
$i=1,\ldots,p-1$. We define a bilinear form on $L_{\mathbb{Z}}$
\begin{eqnarray}\label{bilinear form}
b_L(\alpha,\beta)=-\frac{1}{p}\sum_{\rho\in
\Gamma_p} \rho
(h_L(\alpha,\beta)).
\end{eqnarray}
Note that $b_L$ takes values in $\frac{1}{p}\mathbb{Z}[\omega]$,
so in general $\{L_{\mathbb{Z}},b_L\}$ is not a
$\mathbb{Z}$-lattice. We call it the {\it associated} module
(resp. lattice) of the $\Z[\omega]$-lattice $L$.\\

\textbf{Remark.} Remark. By the definition of the bilinear form is
clear that
$$b_L(\alpha,\beta)=-\frac{1}{p}Tr_{\Q(\omega)/\Q}(h_L(\alpha,\beta)).$$
For a precise definition of the Trace see \cite[page 128]{Ebeling}
\subsection{The $\mathbb{Z}$-lattice $F_p$.}\label{subsection F_p}
We consider a K3 surface admitting an elliptic fibration. Let
$p$ be an odd prime number. Let $I_{p}$ be a semistable fiber
of a minimal elliptic fibration, i.e. (cf. section \ref{elliptic:
general facts}) $I_{p}$ is a fiber which is a reducible curve,
whose irreducible components are the
edges of a $p$-polygon, as described in \cite[Table I.4.1]{miranda1}, we denote the $p$-irreducible components by $C_i$, $i=0,\ldots, p-1$, then
$$
\begin{array}{cc}
C_i\cdot C_j=&\left\{\begin{array}{lll} -2 \ &if\ i\equiv j& \mod p\\
1& if\ |i-j|\equiv 1 &\mod p\\ 0\ & \mbox{otherwise.}&
\end{array}\right.
\end{array}
$$
We consider now the free $\Z$-module $F_p$ with basis the elements of the
form $C_i-C_{i+1}$, $i=1,\ldots,p-1$ and with bilinear form  $b_{F_p}$ which is the restriction of the intersection form to the basis $C_i-C_{i+1}$, then
$\{F_p,b_{F_p}\}$ is a $\Z$-lattice.

\subsection{The $\mathbb{Z}[\omega_p]$-lattice $G_p$.} Let $G_p$ be the
$\mathbb{Z}[\omega]$-lattice
$G_p:=(1-\omega)^2\mathbb{Z}[\omega]$, with the standard hermitian
form: $h(\alpha,\beta)=\alpha\bar{\beta}$. A basis for the
$\mathbb{Z}$-module $G_{p,\mathbb{Z}}$ is
$(1-\omega)^2\omega^i$, $i=0,\ldots,p-2$.\\

On $\mathbb{Z}[\omega]$ we consider the bilinear form $b_L$
defined in \eqref{bilinear form}, with values in
$\frac{1}{p}\mathbb{Z}$,

$$
b(\alpha,\beta)=-\frac{1}{p}\sum_{\rho\in
\Gamma_p}\rho (\alpha
\bar{\beta}),\ \ \ \ \ \ \alpha,\beta\in\mathbb{Z}[\omega],
$$

then we have
\begin{lemma}\label{lemma b_f}
The bilinear form $b$ restricted to $G_p$ (denoted by $b_G$) has values in
$\mathbb{Z}$ and coincides with the intersection form on $F_{p}$
by using the map $F_p\rightarrow G_p$ defined by
$C_i-C_{i+1}\mapsto \omega^i(1-\omega)^2$, $i=1,\ldots, p-1$,
$C_{p}=C_0$.
\end{lemma} \textit{Proof.} An easy computation shows that we have for $p>3$:
$$
\begin{array}{ll}
b_{G}(\omega^k(1-\omega)^2,\omega^h(1-\omega)^2)=&\left\{\begin{array}{ll}
-6\ \ &\mbox{if}\ k\equiv h\ \mod p,\\
4\ \  &\mbox{if}\ |k-h|\equiv 1\ \mod p,\\
-1\ \ &\mbox{if}\ |k-h|\equiv 2\ \mod p,\\
0\ \ &\mbox{otherwise},\\
\end{array}\right.
\end{array}
$$
and for $p=3$:
$$
\begin{array}{ll}
b_{G}(\omega^k(1-\omega)^2,\omega^h(1-\omega)^2)=&\left\{\begin{array}{ll}
-6\ \ &\mbox{if}\ k\equiv h\ \mod p,\\
3\ \ &\mbox{if}\ |k-h|\equiv 1\ \mod p.
\end{array}\right.
\end{array}
$$
The intersection form on $F_p$ is easy to compute (cf. section
\ref{subsection F_p}) and this computation proves that the map
$F_p\rightarrow G_p$ defined in the lemma is an isometry. \eprf


\smallskip

In section \ref{elliptic} we apply the results of this section and we find a
$\mathbb{Z}[\omega]$-lattice $\{L_p,h_{L_p}\}$ such that
\begin{itemize}
\item $\{L_p,h_{L_p}\}$ contains $G_p$ as sublattice; \item
$\{L_{p,\mathbb{Z}},b_{L_p}\}$ is a $\mathbb{Z}$-lattice; \item
the $\mathbb{Z}$-lattice $\{L_{p,\mathbb{Z}},b_{L_p}\}$ is
isometric to the $\mathbb{Z}$-lattice
$(H^2(X,\mathbb{Z})^{\sigma_p^*})^{\perp}$ for $p=3,5,7.$
\end{itemize}

\section{Some general facts on elliptic fibrations}\label{elliptic: general facts}
In the next section we give explicit examples of K3 surfaces
admitting a
symplectic automorphism $\sigma_p$ by using elliptic fibrations.
Here we recall some general results about these fibrations.\\
Let $X$ be an elliptic $K3$ surface, this means that we have a
morphism
\begin{eqnarray*}
f:X\lra \PP^1
\end{eqnarray*}
such that the generic fiber is a (smooth) elliptic curve. We
assume moreover that we have a section $s:\PP^1\lra X$. The
sections of $X$ generate the Mordell-Weil group $MW_f$ of $X$ and
we take $s$ as zero section. This group acts on $X$ by translation
(on each fiber), hence it leaves the two form invariant. We assume
that the singular fibers of the fibration are all of type $I_{m}$,
$m\in\mathbb{N}$. Let $F_j$  be a fiber of type $I_{m_j}$, we
denote by $C_0^{(j)}$ the irreducible components of the fibers
meeting the zero section. After choosing an orientation, we denote
the other irreducible components of the fibers by
$C_1^{(j)},\ldots, C_{m_j-1}^{(j)}$. In the sequel we always
consider $m_j$ a prime number, and the notation $C^{(j)}_i$ means
$i\in\mathbb{Z}/m_j{\Z}$. For each section $r$ we define the
number $k:=k_j(r)$ by
\begin{eqnarray*}
r\cdot C_j^{(k)}=1\ \mbox{and}\ r\cdot C_j^{(i)}=0\ \mbox{if}\
i=0,\ldots, m_j-1\ \ i\neq k.
\end{eqnarray*}
If the section $r$ is a torsion section and $h$ is the number of
reducible fibers of type $I_{m_j}$, then by \cite[Proposition
3.1]{miranda2} we have
\begin{eqnarray}\label{component numbers}
\sum_{j=1}^hk_j(r)\left(1-\frac{k_j(r)}{m_j}\right)=4.
\end{eqnarray}
Moreover we recall the Shioda-Tate formula (cf.\cite[Corollary
5.3]{shioda} or \cite[p.70]{miranda1})
\begin{eqnarray}\label{rank NS}
\mbox{rank}(NS(X))=2+\sum_{j=1}^h(m_j-1)+\mbox{rank}(MW_f).
\end{eqnarray}
The $\mbox{rank}(MW_f)$ is the number of generators of the free
part. If there are no sections of infinite order then
$\mbox{rank}(MW_f)=0$. Assume that $X$ has $h$ fibers of type
$I_m$, $m\in\N$, $m>1$, and the remaining singular fibers are of
type $I_1$, which are rational curves with one node. Let $U\oplus
(A_{m-1})^h$ denote the lattice generated by the zero section, the
generic fiber and by the components of the reducible fibers not
meeting $s$. If there are no sections of infinite order then it
has finite index in $NS(X)$ equal to $n$, the order of the torsion
part of the group $MW_f$. Using this remark we find that
\begin{eqnarray}\label{det NS}
|\det(NS(X))|=\frac{\det(A_{m-1})^h}{n^2}=\frac{m^h}{n^2}.
\end{eqnarray}

\subsection{Elliptic fibrations with a symplectic automorphism.}\label{lenostre}
Now we describe three particular elliptic fibrations which admit a
symplectic automorphism $\sigma_3$, $\sigma_5$ or $\sigma_7$.
Assume that we have a section of prime order $p=3,5,7$. By
\cite[No. 560, 2346, 3256]{shimada} there exist elliptic
fibrations with one of the following configurations of
components of singular fibers $I_p$ not meeting $s$ such that all
the singular fibers of the fibrations are semistable (i.e. they
are all of type $I_n$ for a certain $n\in\mathbb{N}$) and the
order of the torsion subgroup of the Mordell-Weil group
$o(MW_f)=p$ :
\begin{eqnarray}\label{elliptic fibration Shimada}
\begin{array}{lll}
p=3:&6A_2&o(MW_f)=3,\\
p=5:&4A_4&o(MW_f)=5,\\
p=7:&3A_6&o(MW_f)=7.
\end{array}
\end{eqnarray}
We can assume that the remaining singular fibers are of type $I_1$.
Since the sum of the Euler characteristic of the fibers must add
up to $24$, these are six, four, resp. three fibers. Observe that
each section of finite order induces a symplectic automorphism of
the same order which corresponds to a translation by the section
on each fiber, we denote it by $\sigma_p$. The nodes of the $I_1$
fibers are then the fixed points of these automorphisms, whence
$\sigma_p$ permutes the $p$ components of the $I_p$ fibers. 
For these fibrations we have rank $NS(X)=14,~18,~20$ and
dimensions of the moduli spaces six, two and zero, which is one
less then the maximal possible dimension of the moduli space we
have given in the Proposition \ref{theorem moduli}.
\subsubsection{Weierstrass forms.} We compute the Weierstrass form for the elliptic fibration described
in \eqref{elliptic fibration Shimada}. When $X$ is a K3
surface then this form is
\begin{eqnarray}\label{formaw}
y^2=x^3+A(t)x+B(t),~~~~~t\in\PP^1
\end{eqnarray}
or in homogeneous coordinates
\begin{eqnarray}\label{formawh}
x_3x_2^2=x_1^3+A(t)x_1x_3^2+B(t)x_3^3
\end{eqnarray}
where $A(t)$ and $B(t)$ are polynomials of degrees eight and
twelve respectively, $x_3=0$ is the line at infinity and also the
tangent to the
inflectional point $(0:1:0)$.\\
\textit{Fibration with a section of order 3.} In this case the
point of order three must be an inflectional point (cf. \cite[Ex. 5,
p.38]{cassels}), we want to determine
$A(t)$ and $B(t)$ in the equation \eqref{formaw}. We start by imposing to a general line
$y=l(t)x+m(t)$ to be an inflectional tangent so the equation of
the elliptic fibration is
$$
\begin{array}{l}
y^2=x^3+A(t)x+B(t),~~~~~t\in\PP^1,~~~~~\mbox{with}\\
\end{array}
\begin{array}{l}
A(t)=\frac{\textstyle 2l(t)m(t)+l(t)^4}{\textstyle 3},~~~~~
B(t)=\frac{\textstyle m(t)^2-l(t)^6}{\textstyle 3^3}.
\end{array}
$$
Since $A(t)$ and $B(t)$ are of degrees eight and twelve, we
have $\deg l(t)=2$ and $\deg m(t)=6$. The section of order three is
$$t\mapsto \left(\frac{l(t)^2}{3},\frac{l(t)^3}{3}+m(t)\right).$$
The discriminant $\Delta=4A^3+27B^2$  of the fibration is
$$
\Delta=\frac{(5l(t)^3+27m(t))(l(t)^3+3m(t))^3}{27}
$$
hence in general it vanishes to the order three on six values of
$t$ and to the order one on other six values. Since $A$ and $B$ in
general do not vanish on these values, this equation parametrizes
an elliptic fibration with six fibers $I_3$ (so we have six curves
$A_2$ not meeting the zero section) and six fibers $I_1$ (cf.
\cite[Table IV.3.1 pag.41]{miranda1}).

\textit{Fibration with a section of order 5.} In the same way we
can compute the Weierstrass form of the elliptic fibration
described in \eqref{elliptic fibration Shimada} with a section of
order five.\\
In \cite{billing} a geometrical condition for the existence of a point
of order five on an elliptic curve is given.
For fixed $t$ let the cubic curve be in the form \eqref{formawh}
then take two arbitrary distinct lines through $O$ which meet the
cubic in two other distinct points each. Call $1$, $4$ the points
on the first line and $2$, $3$ the points on the second line, then
$1$ (or any of the other point) has order five if:\\
-the tangent through $1$ meets the cubic in $3$,\\
-the tangent through $4$ meets the cubic in $2$,\\
-the tangent through $3$ meets the cubic in $4$,\\
-the tangent through $2$ meets the cubic in $1$.\\
These conditions give the Weierstrass form:
\begin{eqnarray*}
\begin{array}{l}
y^2=x^3+A(t)x+B(t),~~~~~t\in\PP^1,~~~~~\mbox{with}\\
\\
A(t)=\frac{\textstyle (-b(t)^4+b(t)^2a(t)^2-a(t)^4-3a(t)b(t)^3+3a(t)^3b(t))}{\textstyle 3},~~\\
\\
B(t)=\frac{\textstyle (b(t)^2+a(t)^2)(19b(t)^4-34b(t)^2a(t)^2+19a(t)^4+18a(t)b(t)^3-18a(t)^3b(t))}{\textstyle 108}
\end{array}
\end{eqnarray*}
where $\deg a(t)=2$, $\deg b(t)=2$. The section of order five is
$$t\mapsto ((2b(t)^2-a(t)^2)2:3(a(t)+b(t))(a(t)-b(t))^2:6)$$ and the discriminant is
$$\Delta=\frac{1}{16}(b(t)^2-a(t)^2)^5(11(b(t)^2-a(t)^2)+4a(t)b(t)).$$
By a careful analysis of the zeros of the discriminant we can
see that the fibration has four fibers $I_5$ and four fibers $I_1$
(cf.
\cite[Table IV.3.1 pag.41]{miranda1}).\\
\textit{Fibration with a section of order 7.} To find the
Weierstrass form we use also in this case
  the results of \cite{billing}.  We explain briefly the idea to find a set of
points of order seven on an elliptic curve. One takes points
$0,3,4$ and $1,2,4$ on two lines in the plane. Then the
intersections of a line through $3$ different from the lines
$\{3,2\}$, $\{3,4\}$, $\{3,1\}$ with the lines $\{1,0\}$ and
$\{2,0\}$ give two new points $-1$ and $-2$. By using the
conditions that the tangent through $1$ goes through $-2$ and the
tangent through $2$ goes through $3$ one can determine a cubic
having a point of order seven which is e.g. $1$. By using these
conditions one can find the equation, but since the computations
are quite involved, we recall the Weierstrass form given in
\cite[p.195]{tate}
\begin{eqnarray*}
y^2+(1+t-t^2)xy+(t^2-t^3)y=x^3+(t^2-t^3)x^2.
\end{eqnarray*}
By a direct check one sees that the point of order seven is
$(0(t),0(t))$. This elliptic fibration has three fibers $I_7$ and three fibers $I_1$.

\section{Elliptic K3 surfaces with an automorphism of prime order}\label{elliptic}
In this section we prove the main theorem:
\begin{theorem}\label{main}
For any K3 surface $X$ with a symplectic automorphism $\sigma_p$
of order $p=2,3,5,7$ the action on $H^2(X,\Z)$ decomposes in the
following way:
$$
\begin{array}{l}
{\bf p=2}:~~~~~H^2(X,\Z)^{\sigma_2^*}=E_8(-2)\oplus U\oplus U\oplus U,~~(H^2(X,\Z)^{\sigma_2^*})^{\perp}=E_8(-2).\\
\\
{\bf p=3}:~~~~~H^2(X,\Z)^{\sigma_3^*}=U\oplus U(3)\oplus U(3)\oplus A_2\oplus A_2\\
(H^2(X,\Z)^{\sigma_3^*})^{\perp}=\left\{
\begin{array}{ll}
&x_i\equiv x_j~~ \mod(1-\omega_3),\\
(x_1,\ldots,x_6)\in(\mathbb{Z}[\omega_3])^{\oplus 6}\ :&\\
&\sum_{i=1}^6x_i\equiv 0~~ \mod(1-\omega_3)^2
\end{array}\right\}=K_{12}(-2)\\
\\
\mbox{with hermitian form}\
h(\alpha,\beta)=\sum_{i=1}^6(\alpha_i\overline{\beta_i}).\\
\end{array}
$$
$$
\begin{array}{l}
{\bf p=5}:~~~~~H^2(X,\Z)^{\sigma_5^*}=U\oplus U(5)\oplus U(5)\\
(H^2(X,\Z)^{\sigma_5^*})^{\perp}=\left\{
\begin{array}{ll}
&x_1\equiv x_2\equiv 2x_3\equiv 2x_4~~ \mod(1-\omega_5),\\
(x_1,\ldots,x_4)\in(\mathbb{Z}[\omega_5])^{\oplus 4}:&\\
 &(3-\omega_5)(x_1+x_2)+x_3+x_4\equiv 0~~ \mod(1-\omega_5)^2
\end{array}
\right\} \\
\\
\mbox{with hermitian form}\
h(\alpha,\beta)=\sum_{i=1}^2\alpha_i\overline{\beta_i}+\sum_{j=3}^4
f\alpha_j\overline{f\beta_j}\ \mbox{where}\
f=1-(\omega_5^2+\omega_5^3).\\
\\
$$
$$
{\bf p=7}:~~~~~H^2(X,\Z)^{\sigma_7^*}=U(7)\oplus
\left(\begin{array}
        {cc}4&1\\
        1&2\end{array}
        \right)\\
(H^2(X,\Z)^{\sigma_7^*})^{\perp}=\left\{
\begin{array}{ll}
&x_1\equiv x_2\equiv 6x_3~~ \mod(1-\omega_7),\\
(x_1,x_2,x_3)\in(\mathbb{Z}[\omega_7])^{\oplus 3}\ :&\\
&(1+5\omega_7)x_1+3x_2+2x_3\equiv 0~~
 \mod(1-\omega_7)^2
\end{array}
\right\}\\
\\
\mbox{with hermitian form}\
h(\alpha,\beta)=\alpha_1\overline{\beta_1}+
f_1\alpha_2\overline{f_1\beta_2}+f_2\alpha_3\overline{f_2\beta_3}\\
\mbox{where}\
f_1=3+2(\omega_7+\omega_7^6)+(\omega_7^2+\omega_7^5)\ and\
f_2=2+(\omega_7+\omega_7^6).
\end{array}
$$
\end{theorem}
In the case $p=3$, $K_{12}(-2)$ denotes the Coxeter-Todd lattice
with the bilinear form multiplied by $-2$.\\ This theorem gives a
complete description of the invariant sublattice
$H^2(X,\Z)^{\sigma_p^*}$ and its orthogonal complement in
$H^2(X,\Z)$ for the symplectic automorphisms $\sigma_p$ of all
possible prime order $p=2,3,5,7$ acting on a K3 surface. The
results about the order two automorphism is proven by Morrison in
\cite[Theorem
5.7]{morrison}.\\
We describe the lattices of the theorem and their hermitian forms
in the sections from \ref{elliptic3} to \ref{section lattice
omega_7}. The proof is the following: we identify the action of
$\sigma_p^*$ on $H^2(X,\Z)$ in the case of $X$ an elliptic K3
surface, this is done in several propositions in these sections,
then we apply \cite[Theorem 4.7]{Nikulin} which assure the
uniqueness of this action.

\subsection{A section of order three}\label{elliptic3}
Let $X$ be a K3 surface with an elliptic fibration which admits a
section of order three described in  \eqref{elliptic fibration
Shimada} of section \ref{elliptic: general facts}. We recall that
$X$ has six reducible fibres of type $I_3$ and six singular
irreducible fibres of type $I_1$. In the preceding section we have
seen that the rank of the N\'eron-Severi group is $14$.
We determine now $NS(X)$ and $T_X$.\\
Let $t_1$ denote the section of order three and $t_2=t_1+t_1$. Let
$\sigma_3$ be the automorphism of $X$ which corresponds to the
translation by $t_1$. It leaves each fiber invariant and
$\sigma_3^*(s)=t_1,\ \sigma_3^*(t_1)=t_2,\ \sigma_3^*(t_2)=s$.
Denoted by $C_0^{(i)},C_1^{(i)},C_2^{(i)}$ the components of the
$i-th$ reducible fiber ($i=1,\ldots,6$), we can assume that
$C_1^{(i)}\cdot t_1=C_2^{(i)}\cdot t_2=C_0^{(i)}\cdot s=1$.
\begin{prop}\label{NS and T order 3}
A $\mathbb{Z}$-basis for the lattice $NS(X)$ is given by
$$s,t_1,t_2,F,C_1^{(1)},C_2^{(1)},C_1^{(2)},C_2^{(2)},
C_1^{(3)},C_2^{(3)},C_1^{(4)},C_2^{(4)},C_1^{(5)},C_2^{(5)}.$$ Let
$U\oplus A_2^6$ be the lattice  generated by the section, the
fiber and the irreducible components of the six fibers $I_3$ which
do not intersect the zero section $s$. It has index three in the
N\'eron-Severi group of $X$, $NS(X)$. The lattice $NS(X)$ has
discriminant $-3^4$ and its discriminant form is
\begin{eqnarray*} \mathbb{Z}_3(\frac{2}{3})\oplus
\mathbb{Z}_3(\frac{2}{3})\oplus \mathbb{Z}_3(\frac{2}{3})\oplus
\mathbb{Z}_3(-\frac{2}{3}).
\end{eqnarray*}
The transcendental lattice $T_X$ is $$T_X=U\oplus U(3)\oplus
A_2\oplus A_2$$ and has a unique primitive embedding in the
lattice $\Lambda_{K3}$.
\end{prop}
\bprf  It is clear that a $\Q$-basis for $NS(X)$ is given by
$s,F,C_1^{(i)}, C_2^{(i)}$, $i=1,\ldots, 6$. This basis generates
the lattice $U\oplus A_2^6$. It has discriminant $d(U\oplus
A_2^6)=-3^6$. We denote by
\begin{eqnarray*}
\begin{array}{ll}
c_i=2C_1^{(i)}+C_2^{(i)},& C=\sum c_i,\\
d_i=C_1^{(i)}+2C_2^{(i)}, & D=\sum d_i.
\end{array}
\end{eqnarray*}
Since we know that $t_1\in NS(X)$ we can write
\begin{eqnarray*}
t_1=\alpha s+\beta F+\sum \gamma_i C_1^{(i)}+\sum \delta_i
C_2^{(i)},\ \ \ ~~\alpha,\beta,\gamma_i,\delta_i\in\Q.
\end{eqnarray*}
Then by using the fact that $t_1\cdot s=t_1\cdot C_2^{(i)}=0$ and
$t_1\cdot C_1^{(i)}=t_1\cdot F=1$ one obtains that $\alpha=1$,
$\beta=2$ and $\gamma_1=-2/3$, $\delta_1=-1/3$ hence
$\frac{1}{3}C\in NS(X)$. A similar computation with $t_2$ shows
that $\frac{1}{3} D\in NS(X)$. So one obtains that
\begin{eqnarray}\label{section3}
\begin{array}{l}
t_1=s+2F-\frac{1}{3} C\in NS(X),\\
t_2=s+2F-\frac{1}{3} D\in NS(X).
\end{array}
\end{eqnarray}
and so
\begin{eqnarray*}
3(t_2-t_1)=\sum_{i=1}^6 (C_1^{(i)}-C_2^{(i)})=C-D.
\end{eqnarray*}
We consider now the  $\mathbb{Q}$-basis for the N\'eron-Severi
group
$$s,t_1,t_2,F,C_1^{(1)},C_2^{(1)},C_1^{(2)},C_2^{(2)},C_1^{(3)},C_2^{(3)},C_1^{(4)},C_2^{(4)},C_1^{(5)},C_2^{(5)}.$$
By computing the matrix of the intersection form respect to this
basis one finds that the determinant is $-3^4$. By the Shioda-Tate
formula we have $|\det(NS(X))|=3^4$. Hence this is a
$\mathbb{Z}$-basis for the N\'eron-Severi group. We add to the
classes which generate $U\oplus A_2^6$ the classes $t_1$ and
$\sigma_3^*(t_1)=t_2$ given in the formula (\ref{section3}). Since
$d(U\oplus A_2^6)=3^6$ and $d(NS(X))=-3^4$ the index of $U\oplus A_2^6$ in
$NS(X)$ is 3. Observe that this is also a consequence of a general result given at the end of section \ref{elliptic: general facts}.\\
The classes
$$v_i=\frac{C_1^{(i)}-C_2^{(i)}-(C_1^{(5)}-C_2^{(5)})}{3},\ \ \ \ i=1,\ldots,4$$
generate the discriminant group, which is
$NS(X)^{\vee}/NS(X)$$\cong(\Z/3\Z)^{\oplus 4}$.\\
These classes are not orthogonal to each other with respect to the
bilinear form, so we  take
$$
w_1=v_1-v_2,~~~~w_2=v_3-v_4,~~~~w_3=v_1+v_2+v_3+v_4,~~~~w_4=v_1+v_2-(v_3+v_4)
$$
which form an orthogonal basis with respect to the bilinear form
with values in $\Q/\Z$. And it is easy to compute that
$w_1^2=w_2^2=w_3^2=2/3,~~~~w_4^2=-2/3$. The discriminant form of
the lattice $NS(X)$ is then
\begin{eqnarray}
\mathbb{Z}_3(\frac{2}{3})\oplus \mathbb{Z}_3(\frac{2}{3})\oplus
\mathbb{Z}_3(\frac{2}{3})\oplus \mathbb{Z}_3(-\frac{2}{3}).
\end{eqnarray}
The transcendental lattice $T_X$ orthogonal to $NS(X)$ has rank
eight. Since $NS(X)$ has signature $(1,13)$, the transcendental
lattice has signature $(2,6)$. The discriminant form of the
transcendental lattice is the opposite of the discriminant form of
the N\'eron-Severi lattice. So the transcendental lattice has
signature $(2,6)$, discriminant $3^4$, discriminant group
$T_X^{\vee}/T_X\cong (\Z/3\Z)^{\oplus 4}$ and discriminant form
$\mathbb{Z}_3(-\frac{2}{3})\oplus\mathbb{Z}_3(-\frac{2}{3})\oplus\mathbb{Z}(-\frac{2}{3})\oplus\mathbb{Z}_3(\frac{2}{3})$.
By \cite[Cor. 1.13.5]{Nikulin bilinear} we have $T=U\oplus T'$
where $T'$ has rank six, signature $(1,5)$ and $T'$ has
discriminant form as before. These data identify $T'$ uniquely
(\cite[Corollary 1.13.3]{Nikulin bilinear}). Hence it is
isomorphic to $ U(3)\oplus A_2\oplus A_2 $ with generators for the
discriminant form
$$
(e-f)/3,~~~(e+f)/3,~~~(A-B)/3,~~~(A'-B')/3,
$$
where $e,f,A,B,A',B'$ are the usual bases of the lattices.\\
The transcendental lattice
$$T_X=U\oplus U(3)\oplus A_2\oplus A_2$$
has a unique embedding in the lattice $\Lambda_{K3}$ by
\cite[Theorem 1.14.4]{Nikulin bilinear} or \cite[Corollary
2.10]{morrison}.\eprf

\subsubsection{The invariant lattice and its orthogonal
complement}
\begin{prop}\label{invariant and orthogonal part 3}
The invariant sublattice of the N\'eron-Severi group is isometric
to $U(3)$ and it is generated by the classes $F$ and
$s+t_1+t_2$.\\
The invariant sublattice $H^2(X,\mathbb{Z})^{\sigma_3^*}$ is isometric
to $U\oplus U(3)\oplus U(3)\oplus A_2 \oplus A_2$.\\
Its orthogonal complement
$\Omega_3:=(H^2(X,\mathbb{Z})^{\sigma_3^*})^{\perp}$ is the
negative definite twelve dimensional lattice $\{\mathbb{Z}^{12},
M\}$ where $M$ is the bilinear form
$$
\begin{array}{cc}
&\left(
\begin{array}{rrrrrrrrrrrr}
-4&2&-3&-2&0&-2&0&-2&0&-2&0&-2\\
2& -4 &3&1&0&1&0&1&0&1&0&1\\
-3&3&-18&0&0&0&0&0&0&0&3&-9\\
-2&1&0&-6&-3&0&0&0&0&0&0&0\\
0&0&0&-3&-4&3&2&0&0&0&0&0\\
-2&1&0&0&3&-6&-3&0&0&0&0&0\\
0 &0 &0&0&2&-3&-4&3&2&0&0&0\\
-2&1&0&0&0&0&3&-6&-3&0&0&0\\
0 &0&0&0&0&0&2&-3&-4&3&2&0\\
-2&1&0&0&0&0&0&0&3&-6&-3&0\\
0 & 0&3&0&0&0&0&0&2&-3&-4&3\\
-2&1&-9&0&0&0&0&0&0&0&3&-6
\end{array}
\right)
\end{array}
$$
and it is equal to the lattice $(NS(X)^{\sigma_3^*})^{\perp}$.\\
The lattice $\Omega_3$ admits a unique primitive embedding in the
lattice $\Lambda_{K3}$.\\
The discriminant of $\Omega_3$ is $3^6$ and its discriminant form
is $(\mathbb{Z}_3(\frac{2}{3}))^{\oplus 6}$.\\
The isometry $\sigma_3^*$ acts on the discriminant group
$\Omega_3^{\vee}/\Omega_3$ as the identity.
\end{prop}
\bprf It is clear that the isometry $\sigma_3^*$ fixes the classes
$F$ and $s+t_1+t_2$. These generate a lattice $U(3)$ (with basis
$F$ and $F+s+t_1+t_2$).\\ The invariant sublattice
$H^2(X,\Z)^{\sigma_3^*}$ contains $T_X$ and the invariant
sublattice of the N\'eron-Severi group. So
$(H^2(X,\mathbb{Z})^{\sigma_3^*})^{\perp}=(NS(X)^{\sigma_3^*})^{\perp}$,
this lattice has signature $(0,12)$ and by \cite[p. 133]{Nikulin}
the discriminant group is $(\Z/3\Z)^{\oplus 6}$. Hence by
\cite[Theorem 1.14.4]{Nikulin bilinear} there is a unique
primitive embedding of $(H^2(X,\mathbb{Z})^{\sigma_3^*})^{\perp}$
in the K3-lattice. By using the orthogonality conditions one finds
the following basis of  $\Omega_3=(NS(X)^{\sigma_3^*})^{\perp}$:
\begin{eqnarray*}
\begin{array}{l}
b_1=t_2-t_1,\ \ b_2=s-t_2,\ \  b_3=F-3C_2^{(5)},\ \
b_{2(i+1)}=C_1^{(i)}-C_2^{(i)},\
i=1,\ldots,5\\
b_{2j+3}=C_1^{(j)}-C_1^{(j+1)},\ j=1,\ldots,4.
\end{array}
\end{eqnarray*}

An easy computation shows that the Gram matrix of this basis is
exactly the matrix $M$ which indeed has determinant $3^6.$\\
Since $H^2(X,\Z)^{\sigma_3^*}\supseteq T_X\oplus
NS(X)^{\sigma_3^*}=U\oplus U(3)\oplus U(3)\oplus A_2\oplus A_2$
and these lattices have the same rank, to prove that the inclusion
is an equality we compare their discriminants. The lattice
$(H^2(X,\mathbb{Z})^{\sigma_3^*})^{\perp}$ has determinant $3^6$.
So the lattice $(H^2(X,\mathbb{Z})^{\sigma_3^*})$ has determinant
$-3^6$ (because these are primitive sublattices of
$H^2(X,\mathbb{Z})$). The lattice $U\oplus U(3)\oplus U(3)\oplus
A_2\oplus A_2$ has determinant exactly $-3^6$, so
$$H^2(X,\Z)^{\sigma_3^*}= U\oplus
U(3)\oplus U(3)\oplus A_2\oplus A_2.$$ Since
$NS(X)^{\vee}/NS(X)\subset \Omega_3^{\vee}/\Omega_3$ the
generators of the discriminant form of the lattice $\Omega_3$ are
classes $w_1,\ldots,w_6$ with $w_1,\ldots ,w_4$ the classes which
generate the discriminant form of $NS(X)$ (cf. the proof of the
Proposition \ref{NS and T order 3}) and
$$w_5=\frac{1}{3}(b_1+2b_2)=\frac{1}{3}(2s-t_1-t_2)\ \ \ \ w_6=\frac{1}{3}(b_1+2b_2-2b_3)=\frac{1}{3}(2s-t_1-t_2-2F+6C_2^{(5)}).$$
These six classes are orthogonal, with respect to the bilinear
form taking values in $\Q/\Z$, and generate the discriminant form.
Their squares are $w_1^2=w_2^2=w_3^2=w_5^2\equiv\frac{2}{3}\ \mod
2\mathbb{Z}$, $w_4^2=w_6^2\equiv -\frac{2}{3}\ \ \mod
2\mathbb{Z}$. By replacing $w_4$, $w_6$ by $w_4-w_6$, $w_4+w_6$
we obtain the discriminant form
$(\mathbb{Z}_3(\frac{2}{3}))^{\oplus 6}$.\\
By computing the image of $w_i$, $i=1,\ldots, 6$ under
$\sigma_3^*$ one finds that $\sigma_3^*(w_i)-w_i\in\Omega_3$. For
example:
$\sigma_3^*(w_5)-w_5=\frac{1}{3}(2t_1-t_2-s)-\frac{1}{3}(2s-t_1-t_2)=t_1-s$
which is an element of $\Omega_3$ (in fact it is orthogonal to $F$
and to $s+t_1+t_2$). Hence the action of $\sigma_3^*$ is trivial
on $\Omega_3^{\vee}/\Omega_3$ as claimed. \eprf In the next two
subsections we apply the results of section \ref{z[omega]
lattices} about the $\mathbb{Z}[\omega]$-lattices to describe the
lattice $\{\Omega_3,M\}$ and to prove that $\Omega_3$ is
isomorphic to the lattice $K_{12}(-2)$, where $K_{12}$ is the
Coxeter-Todd lattice (cf. e.g. \cite{Coxeter Todd}, \cite{Conway
Sloane} for a description of this lattice).

\subsection{The lattice $\Omega_3$.}\label{section lattice omega_3}

Let $\omega_3$ be a primitive third root of the unity. In this
section we prove the following result (we use the same notations
of section \ref{z[omega] lattices}):
\begin{theorem}
The lattice $\Omega_3$ is isometric to the $\mathbb{Z}$-lattice
associated to the $\mathbb{Z}[\omega_3]$-lattice $\{L_3,\ h_{L_3}\}$ where
$$
L_3=\left\{
\begin{array}{ll}
&x_i\equiv x_j\ \mod(1-\omega_3),\\
(x_1,\ldots,x_6)\in(\mathbb{Z}[\omega_3])^{\oplus 6}\ :&\\
&\sum_{i=1}^6x_i\equiv 0\ \mod(1-\omega_3)^2
\end{array}
\right\}
$$
and $h_{L_3}$ is the restriction of the standard hermitian form on
$\mathbb{Z}[\omega_3]^{\oplus 6}$.
\end{theorem}
\bprf Let $F=F_3^6$ be the $\mathbb{Z}$-sublattice of $NS(X)$ generated
by
$$
C_i^{(j)}-C_{i+1}^{(j)},\ \ i=0,1,2,\ \ \ \ j=1,\ldots,6
$$
with bilinear form induced by the intersection form on $NS(X)$.\\
Let $G=G_3^6$ denote the $\mathbb{Z}[\omega_3]$-lattice
$(1-\omega_3)^2\mathbb{Z}[\omega_3]^{\oplus 6}$ with the standard
hermitian form. This is a sublattice of $\mathbb{Z}[\omega_3]^{\oplus 6}$.
Applying to each component of $G$ the Lemma \ref{lemma b_f} we
know that $\{G_{\mathbb{Z}}, b_G\}$ is a $\mathbb{Z}$-lattice
isometric to the lattice $F$. The explicit isometry is given by
$$
\begin{array}{l}
C_i^{(1)}-C_{i+1}^{(1)}\mapsto
(1-\omega_3)^2(\omega_3^{i-1},0,0,0,0,0)\\
C_i^{(2)}-C_{i+1}^{(2)}\mapsto
(1-\omega_3)^2(0,\omega_3^{i-1},0,0,0,0)\\
\vdots\\
 C_i^{(6)}-C_{i+1}^{(6)}\mapsto
(1-\omega_3)^2(0,0,0,0,0,\omega_3^{i-1}).\\
\end{array}
$$
The multiplication by $\omega_3$ of an element $(1-\omega_3)^2e_j$
(where $e_j$ is the canonical basis) corresponds to a translation
by $t_1$ on a singular fiber, which sends the curve $C_i^{(j)}$ to
the curve $C_{i+1}^{(j)}$. Hence we have a commutative diagram:
$$
\begin{array}{ccccc}
&F&\longrightarrow&G&\\
\sigma_3^*&\downarrow&&\downarrow&\cdot\omega_3\\
&F&\longrightarrow&G.&
\end{array}
$$

The elements $C_i^{(j)}-C_k^{(j)}$, $i,k=0,1,2$, $j=1,\ldots,6$
are all contained in the lattice $\Omega_3=(NS(X)^{\sigma_3^*})^{\perp}$,
but they
do not generate this lattice.
A set of generators for $\Omega_3$ is
$$
s-t_1,\ \ \ t_1-t_2,\ \ \ C_i^{(j)}-C_h^{(k)}\ \ i,h=0,1,2,\
j,k=1,\ldots,6.
$$
From the formula (\ref{section3})
we obtain that
$$
s-t_1=\sum_{j=1}^6[\frac{1}{3}(C_1^{(j)}-C_2^{(j)})+\frac{1}{3}\sigma_3^*(C_1^{(j)}-C_2^{(j)})].
$$
After the identification of $F$ with $G_{\mathbb{Z}}$ we have
$$
s-t_1=(1-\omega_3)^2(\frac{1}{3}(1+\omega_3))(1,1,1,1,1,1)=(1,1,1,1,1,1).
$$
The divisor $t_1-t_2$, which is the image of $s-t_1$ under the
action of $\sigma_3^*$, corresponds to the vector
$(\omega_3,\omega_3,\omega_3,\omega_3,\omega_3,\omega_3).$
Similarly one can see that the element $C_1^{(1)}-C_1^{(2)}$
corresponds to the vector $(1-\omega_3)(1,-1,0,0,0,0,0)$ and more
in general $C_{i}^{(j)}-C_i^{(k)}$ with $j\neq k$ corresponds to
the vector $(1-\omega_3)(\omega_3^{i-1}e_j-\omega_3^{i-1}e_k)$ where $e_i$
is the standard basis. The lattice $L_3$ generated by the vectors
of $G_{\Z}$ and by
$$
\omega_3^i(1,1,1,1,1,1)\ \ \
(1-\omega_3)\omega_3^{i-1}(-e_j+\omega_3e_k)\ \ i=0,1,2,\ \
j,k=1,\ldots 6,
$$
is thus isometric to $\Omega_3$.\\
In conclusion a basis for $L_3$ is
$$
\begin{array}{ll}
l_1=-\omega_3(1,1,1,1,1,1)&
l_2=-\omega_3^2(1,1,1,1,1,1)=(1+\omega_3)(1,1,1,1,1,1)\\
l_3=(1-\omega_3)^2(0,0,0,0,1-\omega_3,0)&
l_4=(1-\omega_3)^2(1,0,0,0,0,0)\\
l_5=(1-\omega_3)(1,-1,0,0,0,0)&
l_6=(1-\omega_3)^2(0,1,0,0,0,0)\\
l_7=(1-\omega_3)(0,1,-1,0,0,0)&
l_8=(1-\omega_3)^2(0,0,1,0,0,0)\\
l_9=(1-\omega_3)(0,0,1,-1,0,0)& l_{10}=(1-\omega_3)^2(0,0,0,1,0,0)\\
l_{11}=(1-\omega_3)(0,0,0,1,-1,0)&
l_{12}=(1-\omega_3)^2(0,0,0,0,1,0).
\end{array}
$$
The identification between $\Omega_3$ and $L_3$ is given by the
map $b_i\mapsto l_i$.\\ After this identification the intersection
form on $\Omega_3$ is exactly the form $b_{|L_3}$ on $L_3$.\\
The basis $l_i$ of $L_3$ satisfies the condition given in the
theorem, and so
$$
\begin{array}{cll}
L_3\subseteq &\{(x_1,\ldots,x_6)\in(\mathbb{Z}[\omega_3])^{\oplus
6}\ :&x_i\equiv x_j\ \mod(1-\omega_3),\\ & &\sum_{i=1}^6x_i\equiv
0\ \mod(1-\omega_3)^2 \}.
\end{array}$$
Since the vectors $(1-\omega_3)^2e_j,$ $(1-\omega_3)(e_i-e_j)$ and
$(1,1,1,1,1,1)$ generate the $\mathbb{Z}[\omega_3]$-lattice
$\{(x_1,\ldots,x_6)\in(\mathbb{Z}[\omega_3])^{\oplus 6}: \ x_i\equiv x_j\ \mod(1-\omega_3),\ \
\sum_{i=1}^6x_i\equiv 0\ \mod(1-\omega_3)^2 \}$ and since they are
all vectors contained in $L_3$, the equality holds.\eprf

\subsection{The Coxeter-Todd lattice $K_{12}$.}

\begin{theorem}
The lattice $\Omega_3$ is
isometric to the lattice $K_{12}(-2)$.
\end{theorem}
\bprf The lattice $K_{12}$ is described by Coxeter and Todd in
\cite{Coxeter Todd} and by Conway and Sloane in \cite{Conway
Sloane}. The lattice $K_{12}$ is the twelve dimensional
$\mathbb{Z}$-module associated to a six dimensional
$\mathbb{Z}[\omega_3]$-lattice $\Lambda_6^{\omega_3}$.\\
The $\mathbb{Z}[\omega_3]$-lattice $\Lambda_6^{\omega_3}$ is
described in \cite{Conway Sloane} in four different ways. We
recall one of them denoted by $\Lambda^{(3)}$ in \cite[Definition
2.3]{Conway Sloane}, which is convenient for us. Let
$\theta=\omega_3-\bar{\omega_3}$, then $\Lambda_6^{\omega_3}$ is
the $\mathbb{Z}[\omega_3]$-lattice
$$\Lambda_6^{\omega_3}=\{(x_1,\ldots,x_6):\ x_i\in\mathbb{Z}[\omega_3], x_i\equiv x_j \mod \theta,\ \sum_{i=1}^6x_i\equiv 0 \mod 3\}$$
with hermitian form $\frac{1}{3} ^tx\bar{y}$. We observe that
$\theta=\omega_3(1-\omega_3)$. The element $\omega_3$ is a unit in
$\mathbb{Z}[\omega_3]$ so the congruence modulo $\theta$ is the
same as the congruence modulo $(1-\omega_3)$. Observing that
$-3=\theta^2$ it is then clear that the
$\mathbb{Z}[\omega_3]$-module $\Lambda_6^{\omega_3}$ is the
$\mathbb{Z}[\omega_3]$-module $L_3$. The $\mathbb{Z}$-modules
$K_{12}$ and $L_{3,{\mathbb{Z}}}$ are isomorphic since they are
the twelve dimensional $\mathbb{Z}$-lattices associated to the
same $\mathbb{Z}[\omega_3]$-lattice. The bilinear form on the
$\mathbb{Z}$-module $K_{12}$ is given by
$$b_{K_{12}}(x,y)=\frac{1}{3}x\bar{y}=\frac{1}{6}Tr(x\bar{y})$$
and the bilinear form on $L_{3,{\mathbb{Z}}}$ is given
by
$$b_{|L_3}(x,y)=-\frac{1}{3}Tr(x\bar{y}).$$
So the $\mathbb{Z}$-lattice $\{L_{3,{\mathbb{Z}}},b_{|L_3}\}$ is
isometric to $K_{12}(-2).$\eprf

\textbf{Remark.} 1) In \cite{Coxeter Todd} Coxeter and Todd give an
explicit basis of the $\Z$-lattice $K_{12}$. By a direct
computation one can find the change of basis between the basis
described in \cite{Coxeter Todd} and the basis $\{b_i\}$ given in the proof
of Proposition \ref{invariant and orthogonal part 3}.\\
2) The lattice $\Omega_3$ does not contain vectors of norm $-2$
(cf. \cite[Lemma 4.2]{Nikulin}), but has 756 vectors of norm $-4$,
4032 of norm $-6$ and 20412 of norm $-8$. Since these properties
define the lattice $K_{12}(-2)$, (cf. \cite[Theorem 1]{Conway
Sloane}), this is another way to prove the equality between
$\Omega_3$ and $K_{12}(-2)$.\\
3) The lattice $K_{12}(-2)$ is generated by vectors of norm $-4$,
\cite[Section 3]{Plesken Pohst}.

\subsection{Section of order five}\label{elliptic5}
Let $X$ be a K3 surface with an elliptic fibration which admits a
section of order five as described in section \ref{elliptic:
general facts}. We recall that $X$ has four reducible fibres of
type $I_5$ and four singular irreducible fibres of type $I_1$. We
have seen that the rank of the N\'eron-Severi group is
$18$. We determine now $NS(X)$ and $T_X$.\\
We label the fibers and their components as described in the
section \ref{elliptic: general facts}. Let $t_1$ denote the
section of order five which meets the first singular fiber in
$C_1^{(1)}$. By the formula (\ref{component numbers}) of section
\ref{elliptic: general facts} up to permutation of the fibers only
the following situations are possible:
\begin{eqnarray*}
t_1\cdot C_1^{(1)}=t_1\cdot C_1^{(2)}=t_1\cdot
C_2^{(3)}=t_1\cdot C_2^{(4)}=1\ \mbox{and}\ \ t_1\cdot C_i^{(j)}=0\ \mbox{otherwise};\\
\mbox{or}\ t_1\cdot C_1^{(1)}=t_1\cdot C_4^{(2)}=t_1\cdot
C_2^{(3)}=t_1\cdot C_3^{(4)}=1\ \mbox{and}\ \ t_1\cdot
C_i^{(j)}=0\ \mbox{otherwise}.
\end{eqnarray*}
Observe that these two cases describe the same situation if we change
the "orientation" on the last two fibers, so we assume to be in
the first case. Let $\sigma_5$ be the automorphism of order five
which leaves each fiber invariant and is the translation by $t_1$,
so $\sigma_5^*(s)=t_1,\ \sigma_5^*(t_1)=t_2,\
\sigma_5^*(t_2)=t_3$, $\sigma_5^*(t_3)=t_4$, $\sigma_5^*(t_4)=s$.

\begin{prop}
A $\mathbb{Z}$-basis for the lattice $NS(X)$ is given by
$$s,t_1,t_2,t_3,t_4,F,C_1^{(1)},C_2^{(1)},
C_3^{(1)},C_4^{(1)},C_1^{(2)},C_2^{(2)},C_3^{(2)},
C_4^{(2)},C_1^{(3)},C_2^{(3)},C_3^{(3)},C_4^{(3)}.$$ Let $U\oplus
A_4^4$ be the lattice  generated by the section, the fiber and the
irreducible components of the four fibers $I_5$ which do not
intersect the zero section $s$. It has index five in the
N\'eron-Severi group
of $X$, $NS(X)$.\\
The lattice $NS(X)$ has discriminant $-5^2$ and its discriminant
form is $$\Z_5(\frac{2}{5})\oplus \Z_5(-\frac{2}{5}).$$ The
transcendental lattice is $$T_X=U\oplus U(5)$$
and has a
unique primitive embedding in the lattice $\Lambda_{K3}.$
\end{prop}
\bprf The proof is similar to the proof of Proposition \ref{NS and
T order 3}. So we sketch it briefly. The classes $s,F,C_i^{(j)}$,
$i=1,\ldots, 4,\ j=1,\ldots,4$, generate $U\oplus A_4^4$. By using
the intersection form, or by the result of \cite[p.
299]{miranda2}, we find
\begin{eqnarray}\label{sections5}
\begin{array}{ll}
t_1=&s+2F-\frac{1}{5}\left[\sum_{i=1}^2
(4C_1^{(i)}+3C_2^{(i)}+2C_3^{(i)}+C_4^{(i)})+\right.\\
 &\left. +\sum_{j=3}^4
(3C_1^{(j)}+6C_2^{(j)}+4C_3^{(j)}+2C_4^{(j)})\right].
\end{array}
\end{eqnarray}
A $\mathbb{Z}$-basis is
$s,t_1,t_2,t_3,t_4,F,C_1^{(1)},C_2^{(1)}$,$C_3^{(1)}$,$C_4^{(1)}$,$C_1^{(2)},C_2^{(2)}$,
$C_3^{(2)}$,
$C_4^{(2)},C_1^{(3)},C_2^{(3)}$,$C_3^{(3)}$,$C_4^{(3)}$. Since
$d(NS(X))=-5^2$ and $d(U\oplus A_4^4)=-5^4$, the index of $U\oplus
A_4^4$ in $NS(X)$ is five.
Let $w_1$ and $w_2$ be
$$
\begin{array}{l}
w_1=\frac{1}{5}(2C^{(1)}_1+4C^{(1)}_2+C^{(1)}_3+3C^{(1)}_4+4C^{(3)}_1+3C^{(3)}_2+2C^{(3)}_3+C^{(3)}_4);\\
w_2=\frac{1}{5}(3C^{(2)}_1+C^{(2)}_2+4C^{(2)}_3+2C^{(2)}_4+C^{(3)}_1+2C^{(3)}_2+3C^{(3)}_3+4C^{(3)}_4).
\end{array}
$$
The classes $v_1=w_1-w_2$, $v_2=w_1+w_2$ are orthogonal classes
and generate the discriminant group of $NS(X)$, the discriminant
form is $$ \Z_5(\frac{2}{5})\oplus \Z_5(-\frac{2}{5}).$$ The
transcendental lattice $T_X$ has rank four, signature $(2,2)$ and
discriminant form $\Z_5(-\frac{2}{5})\oplus \Z_5(\frac{2}{5})$.
Since in this case $T_X$ is uniquely determined by signature and
discriminant form (cf. \cite[Corollary 1.13.3]{Nikulin bilinear})
this is the lattice
$$
T_X=U\oplus U(5).
$$
The transcendental lattice has a unique embedding in the lattice
$\Lambda_{K3}$ by \cite[Theorem 1.14.4]{Nikulin bilinear} or
\cite[Corollary 2.10]{morrison}. \eprf

\subsubsection{The invariant lattice and its orthogonal
complement}
\begin{prop}\label{invariant and orthogonal part 5}
The invariant sublattice of the N\'eron-Severi lattice is
isometric to the lattice $U(5)$ and it is generated by the classes
$F$ and
$s+t_1+t_2+t_3+t_4$.\\
The invariant lattice $H^2(X,\mathbb{Z})^{\sigma_5^*}$ is
isometric to $U\oplus U(5)\oplus U(5)$ and its orthogonal
complement $\Omega_5=(H^2(X,\mathbb{Z})^{\sigma_5^*})^{\perp}$ is
the negative definite sixteen dimensional lattice
$\{\mathbb{Z}^{16}, M\}$ where $M$ is the bilinear form $$
\begin{array}{cc}
&\left(\begin{array}{rrrrrrrrrrrrrrrr}
-4& 2& 0& 0& 0&-1& 0&0&0&-1& 0& 0& -1& 1&-1& 0\\
2&-4& 2& 0& 5& 2&-1& 0& 0& 2&-1& 0& 1&-1& 1& 1\\
0&  2& -4&  2& -5& -1&  2& -1&  0& -1&  2&  -1& 1&-1& 0&-1\\
0&  0&  2& -4&  0&  0& -1&  2&  0&  0& -1&  2&-1& 1& 1&-1\\
0&  5& -5&  0&-50&  0&  0&  0&  0&  0&  0&  0&0& 0& 5  & -15\\
-1&  2& -1&  0&  0& -6&  4& -1& -3&  0&  0& 0&0& 0& 0& 0 \\
0& -1&  2& -1&  0&  4& -6&  4&  1&  0&  0&  0& 0& 0& 0& 0\\
0& 0&-1& 2& 0&-1& 4&-6& 0& 0& 0& 0& 0& 0& 0& 0 \\
0& 0& 0& 0& 0&-3& 1& 0&-4& 3&-1& 0& 2& 0& 0& 0 \\
-1& 2&-1& 0& 0& 0& 0& 0& 3&-6& 4&-1& -3& 0& 0& 0\\
0&-1& 2&-1& 0& 0& 0& 0&-1& 4&-6& 4& 1& 0& 0& 0  \\
0& 0&-1& 2& 0& 0& 0& 0& 0&-1& 4& -6& 0& 0& 0& 0  \\
-1& 1& 1&-1& 0& 0& 0& 0& 2&-3& 1&0& -4& 3&-1& 0   \\
1&-1&-1& 1& 0& 0& 0& 0& 0& 0& 0& 0& 3&-6& 4&-1  \\
-1& 1& 0& 1& 5& 0& 0& 0& 0& 0& 0&0&   -1& 4&-6& 4 \\
0& 1&-1&-1 &  -15& 0& 0& 0& 0& 0& 0& 0&   0&-1& 4&-6 \\
\end{array}
\right)
\end{array}
$$
and it is equal to the lattice $(NS(X)^{\sigma_5^*})^{\perp}$.\\
The lattice $\Omega_5$ admits a unique primitive embedding in the
lattice $\Lambda_{K3}$.\\
The discriminant of $\Omega_5$ is $5^4$ and its discriminant form
is $(\Z_5(\frac{2}{5}))^{\oplus 4}$.\\
The isometry $\sigma_5^*$ acts on the discriminant group
$\Omega_5^{\vee}/\Omega_5$ as the identity.
\end{prop}
\bprf As in the case of an elliptic fibration with a section of
order three, it is clear that $\sigma_5^*$ fixes the classes $F$
and $s+t_1+t_2+t_3+t_4$. These classes generate the lattice
$U(5)$, and so $H^2(X,\mathbb{Z})^{\sigma_5^*}\supseteq U(5)\oplus
T_X=U(5)\oplus U(5)\oplus U$. Using Nikulin's result in \cite[p.
133]{Nikulin} we find that the lattice
$H^2(X,\mathbb{Z})^{\sigma_5^*}$ has determinant $-5^4$, which is
exactly the determinant of $U(5)\oplus U(5)\oplus U$. Since
these have the same rank, we conclude
that $H^2(X,\mathbb{Z})^{\sigma_5^*}=U(5)\oplus U(5)\oplus U$.\\
The orthogonal complement
$(H^2(X,\mathbb{Z})^{\sigma_5^*})^{\perp}$ is equal to
$(NS(X)^{\sigma_5^*})^{\perp}$ as in Proposition \ref{invariant and
orthogonal part 3}. It has signature $(0,16)$ and by \cite[p.
133]{Nikulin} the discriminant group is $(\Z/5\Z)^{\oplus 4}$.
Hence by \cite[Theorem 1.14.4]{Nikulin bilinear} there is a unique
primitive embedding of $(H^2(X,\mathbb{Z})^{\sigma_5^*})^{\perp}$
in the K3-lattice. By using the orthogonality conditions one finds
the following basis of $\Omega_{5}=(NS(X)^{\sigma_5^*})^{\perp}$:
\begin{eqnarray*}
\begin{array}{l}
b_1=s-t_1,\
b_2=t_1-t_2,\ b_3=t_2-t_3,\ b_4=t_3-t_4,\ b_5=F-5C_4^{(3)},\\
b_i=C_{i-5}^{(1)}-C_{i-4}^{(1)},\ i=6,7,8,\ \ \ \
b_{9}=C_1^{(1)}-C_1^{(2)},\\
b_i=C_{i-9}^{(2)}-C_{i-8}^{(2)},\ i=10,11,12,\ \ \ \
b_{13}=C_1^{(2)}-C_1^{(3)},\\
b_i=C_{i-13}^{(3)}-C_{i-12}^{(3)},\ i=14,15,16.
\end{array}
\end{eqnarray*}
The Gram matrix of this basis is exactly the matrix $M$.\\
The generators of the discriminant group of $\Omega_5$ are the
classes $v_1$, $v_2$ of the discriminant form of $NS(X)$ and the
classes
\begin{eqnarray*}
\begin{array}{l}
v_3=\frac{1}{5}(b_3+2b_1+3b_4+4b_2),\\
v_4=\frac{1}{5}(b_3+2b_1+3b_4+4b_2-b_5).
\end{array}
\end{eqnarray*}
These have $v_3^2=-2/5 \mod 2 \Z$, $v_4^2=2/5 \mod 2\Z$. The
generators $v_1$, $2v_2-4v_3-v_4$, $2v_3$, $v_4$ are orthogonal to
each other and have self-intersection $2/5$. \eprf

\subsection{The lattice $\Omega_5$.}\label{section lattice omega_5}
Let $\omega_5$ be a primitive fifth root of the unity. In this
section we prove the following result
\begin{theorem}
The lattice $\Omega_5$ is isometric to the $\mathbb{Z}$-lattice
associated to the $\mathbb{Z}[\omega_5]$-lattice $\{L_5,\
h_{L_5}\}$ where
$$L_5=\left\{
\begin{array}{ll}
&x_1\equiv x_2\equiv 2x_3\equiv 2x_4\ \ \mod(1-\omega_5)\\
(x_1,\ldots,x_4)\in(\mathbb{Z}[\omega_5])^{\oplus 4}:&\\
 &(3-\omega_5)x_1+(3-\omega_5)x_2+x_3+x_4\equiv 0\ \
 \mod(1-\omega_5)^2
\end{array}
\right\} $$
with the hermitian form
\begin{eqnarray}\label{hermitian form 5}h_{L_5}(\alpha,\beta)=
\sum_{i=1}^2\alpha_i\bar{\beta_i}+\sum_{j=3}^4f\alpha_j\overline{f\beta_j}=
\sum_{i=1}^2\alpha_i\bar{\beta_i}+\tau\sum_{j=3}^4\alpha_j\overline{\beta_j}\
,\end{eqnarray} where $\alpha,\beta\in L_5
\subset\mathbb{Z}[\omega_5]^{\oplus 4}$,
$f=1-(\omega_5^2+\omega_5^3)$ and
$\tau=f\overline{f}=2-3(\omega_5^2+\omega_5^3)$.
\end{theorem}
\bprf The strategy of the proof is the same as in the case with an
automorphism of order three, but the situation is more complicated
because the section $t_1$ does not meet all the fibers $I_5$ in
the same component. For this reason the hermitian form of the
$\mathbb{Z}[\omega_5]$-lattice $L_5$ is not the standard hermitian
form on all the components. It is possible to repeat the
construction used in the case of order three, but with the
hermitian form \eqref{hermitian form 5}. We explain now
how we find this hermitian form.\\
Let $F:=F_5^4$ be the lattice generated by the elements
$C_i^{(j)}-C_{i+1}^{(j)}$,  $i=0,\ldots, 4,$ $j=1,\ldots,4$. This
is a sublattice of $(NS(X)^{\sigma_5*})^{\perp}$. A basis is
$$
\begin{array}{lll}
d_{1+i}=(\sigma_5^*)^i(C_1^{(1)}-C_2^{(1)}),&
d_{5+i}=(\sigma_5^*)^i(C_1^{(2)}-C_2^{(2)}),&
d_{9+i}=(\sigma_5^*)^i(C_1^{(3)}-C_2^{(3)}),\\
d_{13+i}=(\sigma_5^*)^i(C_1^{(4)}-C_2^{(4)}),&
 i=0,\ldots,3 &\\
\end{array}
$$
and the bilinear form is thus the diagonal block matrix
$Q=\mbox{diag}(A,A,B,B)$
$$
\begin{array}{cc}
A=\left(\begin{array}{rrrr}-6&4&-1&-1\\4&-6&4&-1\\-1&4&-6&4\\-1&-1&4&-6\end{array}\right),\
&B=\left(\begin{array}{rrrr}-6&-1&4&4\\-1&-6&-1&4\\4&-1&-6&-1\\4&4&-1&-6\end{array}\right).\
\end{array}
$$
We want to identify the multiplication by $\omega_5$ in the
lattice $G$ with the action of the isometry $\sigma_5^*$ on the
lattice $F$. We consider the $\mathbb{Z}[\omega_5]$-module
$G=(1-\omega)^2\mathbb{Z}[\omega]^{\oplus 4}$. Now we consider the
$\mathbb{Z}$-module $G_{\mathbb{Z}}$. The map
$$
\begin{array}{llll}
\phi:&(\sigma_5^*)^i(C_1^{(1)}-C_2^{(1)})&\mapsto&
(1-\omega_5)^2\omega_5^i(1,0,0,0)\\
&(\sigma_5^*)^i(C_1^{(2)}-C_2^{(2)})&\mapsto&(1-\omega_5)^2 \omega_5^i(0,1,0,0)\\
&(\sigma_5^*)^i(C_1^{(3)}-C_2^{(3)})&\mapsto&(1-\omega_5)^2 \omega_5^i(0,0,1,0)\\
&(\sigma_5^*)^i(C_1^{(4)}-C_2^{(4)})&\mapsto&(1-\omega_5)^2 \omega_5^i(0,0,0,1)
\end{array}
$$
is an isomorphism between the $\mathbb{Z}$-modules
$G_{\mathbb{Z}}$ and $F$.\\
Now we have to find a bilinear form $b_G$ on $G$ such that
$\{G_{\mathbb{Z}},b_G\}$ is isometric to $\{F,Q\}$. On the first
and second fiber the action of $\sigma_5^*$ is
$\sigma_5^*(C_i^{(j)})=C_{i+1}^{(j)}$, $j=1,2$, $i=0,\ldots4$, so
$(\sigma_5^*)^i(C_1^{(j)}-C_2^{(j)})=C_{i+1}^{(j)}-C_{i+2}^{(j)}$.
Hence the map $\phi$ operates on the first two fibers in the
following way:
$$
\begin{array}{llll} \phi:&C_{i+1}^{(1)}-C_{i+2}^{(1)}&\mapsto&
(1-\omega_5)^2\omega_5^i(1,0,0,0)\\
&C_{i+1}^{(2)}-C_{i+2}^{(2)}&\mapsto& (1-\omega_5)^2\omega_5^i(0,1,0,0).
\end{array}
$$
This identification is exactly the identification described in
Lemma \ref{lemma b_f}, so on these generators of the lattices $F$
and $G$ we can choose exactly the form described in the lemma.\\
On the third and fourth fiber the action of $\sigma_5^*$ is
different (because $\sigma_5^*$ is the translation by $t_1$ and it
meets the first and second fiber in the component $C_1$ and the
third and fourth fiber in the component $C_2$ ). In fact
$(\sigma_5^*)^i(C_1^{(j)}-C_2^{(j)})=C_{2i+1}^{(j)}-C_{2i+2}^{(j)}$,
$j=3,4$, $i=0,\ldots,4$ and so
\begin{eqnarray}\label{order5 34fiber}\begin{array}{llll}
\phi:&C_{2i+1}^{(3)}-C_{2i+2}^{(3)}&\mapsto&
(1-\omega_5)^2\omega_5^i(0,0,1,0)\\
&C_{2i+1}^{(4)}-C_{2i+2}^{(4)}&\mapsto& (1-\omega_5)^2\omega_5^i(0,0,0,1).
\end{array}\end{eqnarray}
A direct verification shows that the map $\phi$ defines an
isometry between the module generated by $(\sigma_p^*)^i(C_1^{(j)}-C_2^{(j)})$, $i=0,\ldots,4$ and $(1-\omega_5)^2\mathbb{Z}[\omega_5]$, ($j=3,4$) if one considers on
$(1-\omega_5)^2\mathbb{Z}[\omega_5]$ the bilinear form associated to the
hermitian form $h(\alpha,\beta)=\tau\alpha\overline{\beta}$ where
$\tau=(2-3(\omega_5^2+\omega_5^3))$. The real number $\tau$ is the
square of $f=1-(\omega_5^2+\omega_5^3)$, so the hermitian form
above is also
$h(\alpha,\beta)=\tau\alpha\overline{\beta}=f\alpha\overline{f\beta}$.
So now we consider the $\mathbb{Z}[\omega_5]$-lattice
$\mathbb{Z}[\omega_5]^{\oplus 4}$ with the hermitian form $h$ given in
(\ref{hermitian form 5}) and $G$ as a sublattice of
$\{\mathbb{Z}[\omega_5]^{\oplus 4},h\}$. We show that
$L_5=\Omega_5$.
We have to add to the lattice $F$ some classes to obtain the
lattice $\Omega_5$, and so we have to add some vectors to the
lattice $G$ to obtain the lattice $L_5$. It is sufficient to add
to $F$ the classes $s-t_1$, $C_1^{(1)}-C_1^{(2)}$,
$C_1^{(2)}-C_1^{(3)}$, $C_1^{(3)}-C_1^{(4)}$ and their images
under $\sigma_5^*$. These classes correspond to the following
vectors in $\mathbb{Z}[\omega_5]^{\oplus 4}$:
$$
\begin{array}{l}
s-t_1=(1,1,c,c),\\
C_1^{(1)}-C_1^{(2)}=(1-\omega_5)(1,-1,0,0),\\
C_1^{(2)}-C_1^{(3)}=(1-\omega_5)(0,1,-(1+\omega_5^3),0),\\
C_1^{(3)}-C_1^{(4)}=(1-\omega_5)(0,0,(1+\omega_5^3),-(1+\omega_5^3))
\end{array}
$$
where $c=\omega_5(2\omega_5^2-\omega_5+2)$. A basis for the
lattice $L_5$ is then
$$
\begin{array}{ll}
l_1=(1,1,c,c)&
l_2=\omega_5 l_1\\
l_3=\omega_5^2l_1&
l_4=\omega_5^3l_1\\
l_5=(1-\omega_5)^2(0,0,2+4\omega_5+\omega_5^2+3\omega_5^3,0)& l_6=(1-\omega_5)^2(1,0,0,0)\\
l_7=\omega_5l_6&
l_8=\omega_5^2l_6\\
l_9=(1-\omega_5)(1,-1,0,0)& l_{10}=(1-\omega_5)^2(0,1,0,0)\\
l_{11}=\omega_5l_{10}&
l_{12}=\omega_5^2l_{10}\\
l_{13}=(1-\omega_5)(0,1,-(1+\omega_5^3),0) &
l_{14}=(1-\omega_5)^2(0,0,1,0)\\
l_{15}=\omega_5l_{14} & l_{16}=\omega_5^2l_{14}.
\end{array}
$$
The identification between $\Omega_5$ and $L_5$ is given by the
map $b_i\mapsto l_i$. After this identification the intersection
form on $\Omega_5$ is exactly the form $b_{|L_5}$ on $L_5$. \eprf

\textbf{Remark.} 1) We recall that the density of a lattice $L$ of
rank $n$ is $\Delta=V_n/\sqrt{\det L}$ where $V_n$ is the volume
of the $n$ dimensional sphere of radius $r$ (called
\textit{packing radius} of the lattice),
$V_n=r^n\pi^{n/2}/(n/2)!$, $r=\sqrt{\mu}/2$ and $\mu$ is the
minimal norm of a vector of the lattice.\\
The density of $\Omega_{5}$ is
$\Delta=\frac{\pi^8}{8!}\frac{1}{5^2}\approx 0.0094$.\\
2) The lattice $\Omega_5$ does not admit vectors of norm $-2$ and
can be generated by vectors of norm $-4$, and a basis is $b_1$,
$b_2$, $b_3$, $b_4$, $b_5-b_{13}-2b_{14}-3b_{15}-4b_{16}$, $b_6$,
$b_7$, $b_8$, $b_9$, $b_{10}+b_{11}$, $b_{11}+b_{12}$,
$b_{10}+b_{11}+b_{12}$, $b_{13}$, $b_{14}+b_{15}$,
$b_{15}+b_{16}$, $b_{14}+b_{15}+b_{16}$.

\subsection{Section of order seven}\label{elliptic7}
Let $X$ be a K3 surface with an elliptic fibration which admits a
section of order seven as described in section \ref{elliptic:
general facts}. We recall that $X$ has three reducible fibres of
type $I_7$ and three singular irreducible fibres of type $I_1$. We
have seen that the rank of the N\'eron-Severi group is
$20$. We determine now $NS(X)$ and $T_X$.\\
We label the fibers and their components as described in the
section \ref{elliptic: general facts}. Let $t_1$ denote the
section of order seven which meets the first fiber in $C_1^{(1)}$.
Again by the formula (\ref{component numbers}) of section
\ref{elliptic: general facts} we have
\begin{eqnarray*}
t_1\cdot C_1^{(1)}=1,\ t_1\cdot C_2^{(2)}=1,\ t_1\cdot
C_3^{(3)}=1,\ \mbox{and}\ t_1\cdot C_i^{(j)}=0\ \mbox{otherwise}.
\end{eqnarray*}
Let $\sigma_7$ denote the automorphism of order seven which leaves
each fiber invariant and is the translation by $t_1$, so
$\sigma_7^*(s)=t_1$, $\sigma_7^*(t_1)=t_2$, $\sigma_7^*(t_2)=t_3$,
$\sigma_7^*(t_3)=t_4$, $\sigma_7^*(t_4)=t_5$,
$\sigma_7^*(t_5)=t_6$, $\sigma_7^*(t_6)=s$.
The proofs of the next two propositions are very similar to those of the similar propositions in the case of the automorphisms of order three and five, so we omit them.

\begin{prop}
A $\mathbb{Z}$-basis for the lattice $NS(X)$ is given by
$$s,t_1,t_2,t_3,t_4,t_5,t_6,F,C_1^{(1)},
C_2^{(1)},C_3^{(1)},C_4^{(1)},C_5^{(1)}, C_6^{(1)},C_1^{(2)},
C_2^{(2)},C_3^{(2)},C_4^{(2)}, C_5^{(2)},C_6^{(2)}.$$ Let $U\oplus
A_6^3$ be the lattice  generated by the section, the fiber and the
irreducible components of the three fibers $I_7$ which do not
intersect the zero section $s$. It has index seven in the
N\'eron-Severi group
of $X$, $NS(X)$.\\
The lattice $NS(X)$ has discriminant $-7$
and its discriminant form is $\Z_7(-\frac{4}{7})$.\\
The transcendental lattice $T_X$ is the lattice $\{\Z^{\oplus
2},\Upsilon\}$ where
$$ \Upsilon:=\left(\begin{array}
        {cc}4&1\\
        1&2\end{array}
        \right)
$$
and it has a unique primitive embedding
in the lattice $\Lambda_{K3}.$
\end{prop}

\subsubsection{The invariant lattice and its orthogonal complement.}
\begin{prop}\label{invariant and orthogonal part 7}
The invariant sublattice of the N\'eron-Severi lattice is
isometric to the lattice $U(7)$ and it is generated by the classes
$F$ and
$s+t_1+t_2+t_3+t_4+t_5+t_6$.\\
The invariant lattice $H^2(X,\mathbb{Z})^{\sigma_7^*}$ is
isometric to $U(7)\oplus T_X$. Its orthogonal complement
$\Omega_7:=(H^2(X,\mathbb{Z})^{\sigma_7^*})^{\perp}$ is the
negative definite eighteen dimensional lattice $\{\mathbb{Z}^{18},
M\}$ where $M$ is the bilinear form $$
\begin{array}{cc}
&\left(
\begin{array}{rrrrrrrrrrrrrrrrrr}
-4&2&0&0&0&0&0&-1&0&0&0&0&0&1&-1&0&0&0\\
2&-4&2&0&0&0&0&2&-1&0&0&0&0&-1&1&1&-1&0\\
0&2&-4&2&0&0&7&-1&2&-1&0&0&0&0&0&-1&1&1\\
0&0&2&-4&2&0&-7&0&-1&2&-1&0&1&-1&0&0&0& -1\\
0&0&0&2&-4&2&0&0&0&-1&2&-1&-1&1&1&-1&0&0\\
0&0&0&0&2&-4&0&0&0&0&-1&2&-1&0&-1&1&1&-1\\
0&0&7&-7&0&0&-98&0&0&0&0&0&0&0&0&0&7&-21\\
-1&2&-1&0&0&0&0&-6&4&-1&0&0&0&0&0&0&0&0\\
0&-1&2&-1&0&0&0&4&-6&4&-1&0& 0&0&0&0&0&0\\
0&0&-1&2&-1&0&0&-1&4&-6&4&-1&0&0&0&0&0&0\\
0&0&0&-1&2&-1&0&0&-1&4&-6&4&-1&0&0&0&0&0\\
0&0&0&0&-1&2&0&0&0&-1&4&-6&3&0&0&0&0&0\\
0&0&0&1&-1&-1&0&0&0&0&-1&3&-4&3&-1&0&0&0\\
1&-1&0&-1&1&0&0&0&0&0&0&0&3&-6&4&-1&0&0\\
-1&1&0&0&1&-1&0&0&0&0&0&0&-1&4&-6&4&-1&0\\
0&1&-1&0&-1&1&0&0&0&0&0&0&0&-1&4&-6&4&-1\\
0&-1&1&0&0&1&7&0&0&0&0&0&0&0&-1&4&-6&4\\
0&0&1&-1&0&-1&-21&0&0&0&0&0&0&0&0&-1&4&-6
\end{array}
\right)
\end{array}
$$
and it is equal to the lattice $(NS(X)^{\sigma_7^*})^{\perp}$.\\
The lattice $\Omega_7$ admits a unique primitive embedding in the
lattice $\Lambda_{K3}$.\\
The discriminant of $\Omega_7$ is $7^3$ and its discriminant form
is $(\Z_7(\frac{4}{7}))^{\oplus 3}$.\\
The isometry $\sigma_7^*$ acts on the discriminant group
$\Omega_7^{\vee}/\Omega_7$ as the identity.
\end{prop}
The basis of $(NS(X)^{\sigma_7^*})^{\perp}$ associated to the
matrix $M$ is $b_1=s-t_1$, $b_2=t_1-t_2$, $b_3=t_2-t_3$,
$b_4=t_3-t_4$, $b_5=t_4-t_5$, $b_6=t_5-t_6$ $b_7=F-7C_6^{(2)}$,
$b_i=C_{i-7}^{(1)}-C_{i-6}^{(1)}$, $i=8,\ldots, 12$,
$b_{13}=C_1^{(1)}-C_1^{(2)}$, $b_i=C_{i-13}^{(2)}-C_{i-8}^{(2)}$,
$i=14,\ldots,18$.

\subsection{The lattice $\Omega_7$.}\label{section lattice omega_7}

Let $\omega_7$ be a primitive seventh root of the unity. In this
section we prove the following result
\begin{theorem}
The lattice $\Omega_7$ is isometric to the $\mathbb{Z}$-lattice
associated  to
the $\mathbb{Z}[\omega_7]$-lattice $\{L_7,\ h_{L_7}\}$ where
$$L_7=\left\{
\begin{array}{ll}
&x_1\equiv x_2\equiv 6x_3\ \ \mod(1-\omega_7),\\
(x_1,x_2,x_3)\in(\mathbb{Z}[\omega_7])^{\oplus 3}\ :&\\
&(1+5\omega_7)x_1+3x_2+2x_3\equiv 0\ \
 \mod(1-\omega_7)^2
\end{array}
\right\}
$$
with the hermitian form
\begin{eqnarray}\label{hermitian form 7}
h_{L_7}(\alpha,\beta)=\alpha_1\bar{\beta_1}+f_1\alpha_2\overline{f_1\beta_2}+f_2\alpha_3\overline{f_2\beta_3},
\end{eqnarray} where
$f_1=3+2(\omega_7+\omega_7^6)+(\omega_7^2+\omega_7^5)$,
$f_2=2+(\omega_7+\omega_7^6)$.
\end{theorem}
\bprf As in the previous cases we define the lattice $F:=F_7^3$.
We consider the hermitian form
\begin{eqnarray*}
h(\alpha,\beta)=\alpha_1\bar{\beta_1}+f_1\alpha_2\overline{f_1\beta_2}+f_2\alpha_3\overline{f_2\beta_3}\end{eqnarray*}
on the lattice $\mathbb{Z}[\omega_7]^{\oplus 3}$, and define $G$
to be the
sublattice $G=(1-\omega_7)^2\mathbb{Z}[\omega_7]^{\oplus 3}$ of $\{\mathbb{Z}[\omega_7]^{\oplus 3},h\}.$\\
The map $\phi: F\rightarrow G$
$$
\begin{array}{llll}
\phi:&(\sigma_7^*)^i(C_1^{(1)}-C_2^{(1)})&\mapsto&(1-\omega_7)^2
\omega_7^i(1,0,0)\\
&(\sigma_7^*)^i(C_1^{(2)}-C_2^{(2)})&\mapsto&(1-\omega_7)^2 \omega_7^i(0,1,0)\\
&(\sigma_7^*)^i(C_1^{(3)}-C_2^{(3)})&\mapsto& (1-\omega_7)^2\omega_7^i(0,0,1)\\
\end{array}
$$
is an isomorphism between the $\mathbb{Z}$-lattice
$G_{\mathbb{Z}}$, with the bilinear form induced by the hermitian
form, and $F$ with the intersection form.
We have to add to $G$ some vectors to find a lattice $L_7$
isomorphic to $\Omega_7$. These vectors are
$$\begin{array}{l}
s-t_1=(1,c,k),\\
C_1^{(1)}-C_1^{(2)}=(1-\omega_7)(1,-(1+\omega_7^4),0),\\
C_1^{(2)}-C_1^{(3)}=(1-\omega_7)(0,(1+\omega_7^4),-(1+\omega_7^3+\omega_7^5)),
\end{array}$$
where $c=1+3\omega_{7}+3\omega_{7}^4+\omega_{7}^5$ and
$k=-5+\omega_{7}-5\omega_{7}^2-3\omega_{7}^4-3\omega_{7}^5$. A
basis for the lattice $L_7$ is
$$
\begin{array}{ll}
l_1=(1,c,k)& l_2=\omega_7l_1\\
l_3=\omega_7^2l_1& l_4=\omega_7^3l_1\\
l_5=\omega_7^4l_1&
l_6=\omega_7^5l_1\\
l_7=(1-\omega_7)^2(0,2+4\omega_7+6\omega_7^2+\omega_7^3+3\omega_7^4+5\omega_7^5,0)&
l_8=(1-\omega_7)^2(1,0,0)\\ l_9=\omega_7l_8&
l_{10}=\omega_7^2l_8\\
l_{11}=\omega_7^3l_8&
l_{12}=\omega_7^4l_8\\
l_{13}=(1-\omega_7)(1,-(1+\omega_7^4),0)&
l_{14}=(1-\omega_7)^2(0,1,0)\\
l_{15}=\omega_7l_{14} &
l_{16}=\omega_7^2l_{14}\\
l_{17}=\omega_7^3l_{14} & l_{18}=\omega_7^4l_{14}
\end{array}.
$$
The identification between $\Omega_7$ and $L_7$ is given by the
map $b_i\mapsto l_i$. After this identification the intersection
form on $\Omega_7$ is exactly the form $b_{|L_7}$ on $L_7$ induced
by the hermitian form (\ref{hermitian form 7}). \eprf
\textbf{Remark.} 1) The density of $\Omega_{7}$ is
$\Delta=\frac{\pi^9}{9!}\frac{1}{\sqrt{7^3}}\approx 0.0044$.\\
2) As in the previous cases the lattice $\Omega_7$ does not admit
vectors of norm $-2$ and can be generated by vectors of norm $-4$,
and a basis is $b_1$, $b_2$, $b_3$, $b_4$, $b_5$, $b_6$,
$b_7-b_{13}-2b_{14}-3b_{15}-4b_{16}-5b_{17}-6b_{18}$, $b_8+b_9$,
$b_9+b_{10}$, $b_{10}+b_{11}$, $b_{11}+b_{12}$,
$b_{10}+b_{11}+b_{12}$, $b_{13}$, $b_{14}+b_{15}$,
$b_{15}+b_{16}$, $b_{16}+b_{17}$, $b_{17}+b_{18}$,
$b_{16}+b_{17}+b_{18}$.

\section{Families of K3 surfaces with a symplectic automorphism of order $p$}\label{minpic}
In the previous sections we used elliptic K3 surfaces to describe
some properties of the automorphism $\sigma_p$. All these K3
surfaces have Picard number $\rho_p+1$, where $\rho_p$ is the
minimal Picard number found in the Proposition \ref{theorem
moduli}. In this section we want to describe algebraic K3 surfaces with
symplectic automorphism of order $p$ and with the minimal possible
Picard number. Recall that the values of $\rho_p$ are
\begin{eqnarray*}
\begin{array}{cccc}
p&3&5&7\\
\rho_p&13&17&19,
\end{array}
\end{eqnarray*}
and $\Omega_p$ denote the lattices described in the sections
\ref{elliptic3}, \ref{elliptic5}, \ref{elliptic7}.
\begin{prop}\label{latticeminimalrank}
Let $X$ be a K3 surface with symplectic automorphism of order
$p=3,5,7$ and Picard number $\rho_p$ as above. Let $L$ be a
generator of $\Omega_p^{\perp}\subset NS(X)$, with $L^2=2d>0$ and
let
\begin{eqnarray*}
\mathcal{L}^p_{2d}:=\Z L\oplus \Omega_p.
\end{eqnarray*}
Then we may assume that $L$ is ample and\\
(1) if $L^2\equiv 2,4,\ldots,2(p-1)$ $\mod 2p$, then $\cl^p_{2d}=NS(X)$,\\
(2) if $L^2\equiv 0$ $\mod 2p$, then either $\cl^p_{2d}=NS(X)$ or
$NS(X)=\widetilde{\cl^p_{2d}}$ with
$\widetilde{\cl^p_{2d}}/\cl^p_{2d}\simeq\Z/p\Z$ and in particular
$\widetilde{\cl^p_{2d}}$ is generated by an element $(L/p,v/p)$
with $v^2\equiv 0$ mod $2p$ and $L^2+v^2\equiv 0$ mod $2p^2$.
\end{prop}
\bprf Since $L^2>0$ by Riemann Roch theorem we can assume $L$ or
$-L$ effective. Hence we assume $L$ effective. Let $N$ be an
effective $(-2)$ curve then $N=\alpha L+v'$, with $v'\in\Omega_p$
and $\alpha>0$ since $\Omega_p$ do not contains $(-2)$-curves. We
have $L\cdot N=\alpha L^2>0$, and so $L$ is ample. Moreover recall
that $L$ and $\Omega_p$ are primitive sublattices of $NS(X)$.
Since the discriminant group of $\mathcal{L}^p_{2d}:=\Z L\oplus
\Omega_p$ is $(\Z/2d\Z)\oplus(\Z/p\Z)^{\oplus n_p}$, with $n_3=6$,
$n_5=4$, $n_7=3$ an element in $NS(X)$ not in $\cl^p_{2d}$ is of
the form $(\alpha L/2d,v/p)$, $v\in\Omega_p$
and satisfy the following conditions:\\
(a) $p\cdot (\alpha L/2d,v/p)\in NS(X)$,\\
(b) $(\alpha L/2d,v/p)\cdot L\in \Z$,\\
(c) $(\alpha L/2d,v/p)^2\in\Z$.\\
By using the condition (a) we obtain $p\cdot (\alpha
L/2d,v/p)-(0,v)\in NS(X)$ and so
$$
\frac{p\alpha L}{2d}\in NS(X).
$$
Hence by the primitivity of $L$ in $NS(X)$ follows that $d\equiv
0$ mod $p$, $d=pd'$, $d'\in\Z_{>0}$ and so
$$
\frac{\alpha L}{2d'}\in NS(X)
$$
which gives $\alpha=2d'$ and the class (if there is) is
$(L/p,v/p)$. Now condition (b) gives
$$
(L/p,v/p)\cdot L=L^2/p\in\Z
$$ and so $L^2=2p\cdot r$, $r\in\Z_{>0}$, since the lattice is even. And so if $NS(X)=\widetilde{\cl^p_{2d}}$,
then $L^2\equiv 0$ mod $2p$. Finally condition (c) gives
$$
(L/p,v/p)^2=\frac{L^2+v^2}{p^2}
$$
and so since a square is even $L^2+v^2\equiv 0$ mod $2p^2$.\eprf
In the sections \ref{elliptic3}, \ref{elliptic5}, \ref{elliptic7}
we defined a symplectic automorphism $\sigma_p$, $p=3,5,7$ of
order $p$ on some special $K3$ surfaces and we found the lattices
$\Omega_p=(\Lambda_{K3}^{\sigma_p^*})^{\perp}$. Now we consider
more in general an isometry on $\Lambda_{K3}$ defined as
$\sigma_p^*$ (we call it again $\sigma_p^*$). In the next theorem
we prove that if $X$ is a K3 surface such that
$NS(X)=\mathcal{L}^p_{2d}$ or $NS(X)=\widetilde{\cl^p_{2d}}$, then
this isometry is induced by a symplectic automorphism of the
surface $X$.
\begin{prop}\label{periodi}
Let $\cl_p=\mathcal{L}^p_{2d}$ or $\cl_p=\widetilde{\cl^p_{2d}}$
if $p=3,5$ and let $\cl_7=\widetilde{\cl^7_{2d}}$. Then there
exists a K3 surface $X$ with symplectic automorphism $\sigma_p$ of
order $p$ such that $NS(X)=\cl_p$ ($p=3,5,7$) and
$(H^2(X,\Z)^{\sigma_p^*})^{\perp}=\Omega_p$.\\ Moreover there are
no K3 surfaces with N\'eron-Severi group isometric to
$\mathcal{L}^7_{14d}$.
\end{prop}
\bprf  Let $\sigma_p^*$, $p=3,5,7$, be an isometry as in the sections
\ref{elliptic3}, \ref{elliptic5}, \ref{elliptic7}.
We make the proof in several steps.\\
\textit{Step 1: there exists a marked K3 surface $X$ such that
$NS(X)$ is isometric to $\cl_p$, and there are no K3 surfaces with
N\'eron-Severi group isometric to $\mathcal{L}^7_{14d}$.} By
\cite[Theorem 1.14.4]{Nikulin bilinear} the lattices $\cl^p_{2d}$
$\widetilde{\cl^3_{2d}}$, $\widetilde{\cl^5_{2d}}$ have a unique
primitive embedding in the K3 lattice. The lattice $T_5=U(5)\oplus
U(5)\oplus\langle -2d\rangle$ has a unique primitive embedding in
$\Lambda_{K3}$, again by \cite[Theorem 1.14.4]{Nikulin bilinear}.
Its signature is $(2,3)$ and its discriminant form is the opposite
of the discriminant form of $\cl_{2d}^5$. Since, by
\cite[Corollary 1.13.3]{Nikulin bilinear}, $\cl_{2d}^5$ is
uniquely determined by its signature and discriminant form, it is
the orthogonal of $T_5$ in $\Lambda_{K3}$ and then $\cl_{2d}^5$
admits a primitive embedding in $\Lambda_{K3}$. The lattice
$\widetilde{\cl_{2d}^7}$ is a primitive sublattice of the
N\'eron-Severi group of the K3 surface described in the section
\ref{elliptic7}, so it is a primitive sublattice of $\Lambda_{K3}$
(the same argument can be applied to the lattices
$\widetilde{\cl_{2d}^p}$, $p=3,5$). Let now
$\omega\in\cl_p^{\perp}\otimes\mathbb{C}\subseteq\Lambda_{K3}\otimes\mathbb{C}$,
with $\omega\omega=0$, $\omega\bar{\omega}>0$. We choose $\omega$
generic with these properties. By the surjectivity of the period
map of K3 surfaces, $\omega$ is the period of a K3 surface $X$
with
$NS(X)=\omega^{\perp}\cap\Lambda_{K3}=\cl_p$.\\
The rank of the lattice $\cl_{14d}^7$ is 19 and its discriminant
group has four generators. If $\cl_{14d}^7$ was the N\'eron-Severi
group of a K3 surface, the transcendental lattice of this surface
should be a rank three lattice with a discriminant group generated
by four elements. This is clearly impossible.\\
\textit{Step 2: the isometry $\sigma_p^*$ fixes the sublattice
$\cl_p$.} Since $\sigma_p^*(\Omega_p)=\Omega_p$ and
$\sigma_p^*(L)=L$ (because $L\in\Omega_p^{\perp}$ which is the
invariant sublattice of $\Lambda_{K3}$), if
$\cl_p=\mathcal{L}^p_{2d}=\Z L\oplus \Omega_p$ it is clear that
$\sigma_p^*(\cl_p)=\cl_p$. Now we consider the case
$\cl_p=\widetilde{\cl^p_{2d}}$. The isometry $\sigma_p^*$ acts
trivially on $\Omega_p^{\vee}/\Omega_p$ (cf. Propositions
\ref{invariant and orthogonal part 3}, \ref{invariant and
orthogonal part 5}, \ref{invariant and orthogonal part 7}) and on
$(\mathbb{Z}L)^{\vee}/\mathbb{Z}L$. Let $\frac{1}{p}(L,v')\in
\cl_p$, with $v'\in\Omega_p$. This is also an element in
$(\Omega_p\oplus L\mathbb{Z})^{\vee}/(\Omega_p\oplus
L\mathbb{Z})$. So we have
$\sigma_p^*(\frac{1}{p}(L,v'))\equiv\frac{1}{p}(L,v')\ \
\mod(\Omega_p\oplus\mathbb{Z}L)$, which means
$$\sigma_p^*(\frac{1}{p}(L,v'))=\frac{1}{p}(L,v')+(\beta L,v''),\qquad \beta \in\mathbb{Z},\ \ v''\in\Omega_p.$$
Hence we have $\sigma_p^*(\cl_p)=\cl_p$.\\
\textit{Step 3: The isometry $\sigma_p^*$ is induced by an
automorphism of the surface $X$.} The isometry $\sigma_p^*$ fixes
the sublattice $\cl_p^{\perp}$ of $\Lambda_{K3}$, so it is an
Hodge isometry. By the Torelli theorem an effective Hodge isometry
of the lattice $\Lambda_{K3}$ is induced by an automorphism of the
K3 surface (cf. \cite[Theorem 11.1]{bpv}). To apply this theorem
we have to prove that $\sigma_p^*$ is an effective isometry. An
effective isometry on a surface $X$ is an isometry which preserves
the set of effective divisors. By \cite[Corollary 3.11]{bpv}
$\sigma_p^*$ preserves the set of the effective divisors if and
only if it preserves the ample cone. So if $\sigma_p^*$ preserves
the ample cone it is induced by an automorphism of the surface.
This automorphism is symplectic by construction (it is the
identity on the transcendental lattice
$T_X\subset\Omega_p^{\perp}$), and so if $\sigma_p^*$ preserves
the ample cone, the theorem is proven.\\
\textit{Step 4: The isometry $\sigma_p^*$ preserves the ample cone
$\mathcal{A}_X$.} Let $\mathcal{C}_X^+$ be one of the two
connected components of the set $\{x\in H^{1,1}(X,\mathbb{R})\ |\
(x,x)>0\}$. The ample cone of a K3 surface $X$ can be described as
the set $\mathcal{A}_X=\{x\in\mathcal{C}_X^+\ |\ (x,d)>0\
\textrm{for each}\ d\ \textrm{such that}\ (d,d)=-2,\ d\
\textrm{effective}\}.$ First we prove that $\sigma_p^*$ fixes the
set of the effective $(-2)$-curves. Since there are no
$(-2)$-curves in $\Omega_p$, if $N\in\cl_p$ has $N^2=-2$ then
$N=\frac{1}{p}(aL,v')$, $v'\in\Omega_p$, for an integer $a\neq 0$.
Since $\frac{1}{p} aL^2=L\cdot N>0$, because $L$ and $N$ are
effective divisor, we obtain $a>0$. The curve $N'=\sigma_p^*(N)$
is a $(-2)$-curve because $\sigma_p^*$ is an isometry, hence $N'$
or $-N'$ is effective. Since
$N'=\sigma_p^*(N)=(aL,\sigma_p^*(v'))$ we have $-N'\cdot
L=-aL^2<0$ and so $-N'$ is not effective. Using the fact that
$\sigma_p^*$ has finite order it is clear that
$\sigma_p^*$ fixes the set of the effective $(-2)$-curves.\\
Now let $x\in\mathcal{A}_X$ then $\sigma_p^*(x)\in\mathcal{A}_X$,
in fact $(\sigma_p^*(x),\sigma_p^*(x))=(x,x)>0$ and for each effective
$(-2)$-curve $d$ there exists an effective $(-2)$-curve $d'$ with
$d=\sigma_p^*(d')$, so we have
$(\sigma_p^*(x),d)=(\sigma_p^*(x),\sigma_p^*(d'))=(x,d')>0$. Hence
$\sigma^*_p$ preserves $\mathcal{A}_X$ as claimed. \eprf
\begin{cor}
The coarse moduli space of $\mathcal{L}_p$-polarized $K3$ surfaces (cf. \cite{dolgachev} for the definition) $p=3,5,7$  has dimension seven, three, respectively one and is a quotient of
\begin{eqnarray*}
\mathcal{D}_{\cl_p}=\{\omega\in\PP(\cl_p^{\perp}\otimes_{\Z}\C):~\omega^2=0,~\omega\bar{\omega}>0\}
\end{eqnarray*}
by an arithmetic group $O(\cl_p)$.
\end{cor}
{\bf Remark.} In particular the moduli space of K3 surfaces admitting a symplectic automorphism of order $p=3,5,7$ has dimension respectively seven, three and one.

\section{Final remarks}
1. In Proposition \ref{periodi} it would be interesting to prove
the unicity of the lattices $\widetilde{\mathcal{L}_{2d}^{p}}$,
this requires some careful analysis of the automorphism group of
the lattices $\Omega_p$, $p=3,5,7$.\\
2. It is not difficult to give examples of K3 surfaces (not elliptic) in some projective space with a symplectic automorphism of order three or five. Consider for example the surfaces of $\mathbb{P}^3$:
\begin{eqnarray*}
\begin{array}{ll}
S1:&q_4(x_0,x_1)+q_2(x_0,x_1)x_2x_3+l_1(x_0,x_1)x_2^3+l_1'(x_0,x_1)x_3^3+ax_2^2x_3^2=0\\
S2:& a_{01}x_0^2x_1^2+a_{23}x_2^2x_3^2+a_{0123}x_0x_1x_2x_3+a_{02}x_0^3x_2+a_{13}x_1^3x_3+a_{12}x_1x_2^3+a_{03}x_0x_3^3=0.\\
\end{array}
\end{eqnarray*}
where $q_i$ is homogeneous of degree $i$, $l_1$, $l_1'$ are linear
forms, and $a_{ij}\in\C$.
The surfaces $S_1$ resp. $S_2$ admit symplectic automorphisms of order three  resp. of order five induced by the automorphisms of $\mathbb{P}^3$ given by $\sigma_3:(x_0:x_1:x_2:x_3)\longrightarrow (x_0:x_1:\omega_3 x_2:\omega_3^2 x_3)$ and $\sigma_5: (x_0:x_1:x_2:x_3)\longrightarrow (\omega_5 x_0:\omega_5^4 x_1:\omega_5^2 x_2:\omega_5^3 x_3)$. The automorphisms of $\PP^3$ commuting with $\sigma_3$ resp. $\sigma_5$ form a space of dimension six, resp. four, since the equations depend on 13, resp. seven parameters the dimension of the moduli space is seven, resp. three as expected (this is the minimal possible dimension). In a similar way one can costruct many more examples. In the case of order seven automorphisms it is more difficult to give such examples. Already in the case of a polarization $L^2=2$, the K3 surface is the minimal resolution of the double covering of $\mathbb{P}_2$ ramified on a sextic with singular points and these are the fixed points of the automorphisms. One should resolve the singularities and analyze the action on the resolution before doing the double cover.

\end{document}